\documentclass[11pt]{amsart}

\usepackage{amsmath}
\usepackage{amssymb}
\usepackage{amsthm}
\usepackage{graphicx} 
\usepackage{url}
\usepackage[bookmarks]{hyperref}
\usepackage[msc-links]{amsrefs} 
\usepackage{multirow}
\usepackage{array}
\usepackage{wasysym}
\usepackage[shortlabels]{enumitem}
\usepackage{color}
\usepackage[normalem]{ulem}
\usepackage{soul}
\usepackage{rotate}
\usepackage{comment}

\usepackage[T1]{fontenc}
\usepackage{fourier}

\theoremstyle{plain}
\newtheorem{introtheorem}{Theorem}

\theoremstyle{plain}
\newtheorem{theorem}{Theorem}[section]

\newtheorem{lemma}[theorem]{Lemma}
\newtheorem{corollary}[theorem]{Corollary}
\newtheorem{problem}[theorem]{Problem}

\theoremstyle{definition}
\newtheorem{definition}[theorem]{Definition}

\newtheorem{question}[theorem]{Question}
\newtheorem{remark}[theorem]{Remark}

\newcommand{\bsm}{\begin{smallmatrix}}
\newcommand{\esm}{\end{smallmatrix}}

\newcommand{\dgood}{\textrm{dist}_{\mathbf{G}}}
\newcommand{\ggood}{\mathbf{G}}
\newcommand{\GLZ}{\mathrm{GL}_2(\ZZ)}
\newcommand{\SLZ}{\mathrm{SL}_2(\ZZ)}

\newcommand{\Phisus}{\Phi_{\mathrm{sus}}}
\newcommand{\PP}{\mathbb{P}}
\newcommand{\RR}{\mathbb{R}}
\newcommand{\TT}{\mathbb{T}}
\newcommand{\ZZ}{\mathbb{Z}}
\newcommand{\tr}{\mathrm{tr}}
\newcommand{\Th}{\TT^2_{2 /3}}
\newcommand{\scup}{\#}

\newcommand{\bpm}{\begin{pmatrix}}
\newcommand{\epm}{\end{pmatrix}}

\newcommand{\term}[1]{{\bf #1}}

\title{Almost equivalence of suspension Anosov flows}

\author{Pierre Dehornoy}
\address{Université Aix-Marseille, CNRS, I2M, 13000 Marseille, France}
\email{pierre.dehornoy@univ-amu.fr}

\author{Mario Shannon}
\address{Facultad de Ingenier\'ia, Universidad de la Rep\'ublica, 11200 Montevideo, Uruguay}
\email{shannon@cmat.edu.uy}
\date{July 2026}

\begin{document}

\begin{abstract}
We provide a written proof of a result due to H. Minakawa, which states that all suspension Anosov flows generated by hyperbolic matrices with positive trace are pairwise almost equivalent. The proof relies on constructing, for any given suspension flow, a genus-one Birkhoff section whose first-return map has fewer fixed points than the original map. We improve Minakawa's result by explicitly calculating the first return map onto this section, which leads to explicit bounds on the distances between suspension Anosov flows within the graph of Anosov flows. 
\end{abstract}

\thanks{We thank Thomas Barthelmé and Chi Cheuk Tsang for some remarks concerning questions and research directions, and the anonymous referee for useful comments. PD was and is supported by the French ANR projects Groméov ANR-19-CE40-0007 and AnoDyn ANR-24-CE40-5065.}

\maketitle

\section{Introduction}
\label{S:Introduction}

A flow generated by a non-singular vector field of class $C^1$ on a closed Riemannian 3-manifold is called \term{Anosov} if its induced action on the tangent bundle of the manifold preserves a \emph{hyperbolic splitting}. This family of dynamical systems has been widely studied since the pioneering works of Hadamard in the 19th century, and has played a central role in the development of \emph{chaos theory}. 
We refer to \cite{Katok-Hasselblatt} for a general account on dynamical systems preserving a hyperbolic splitting. 

There are two basic constructions of such flows: \emph{suspensions of hyperbolic automorphisms of the 2-dimensional torus} and \emph{geodesic flows on closed hyperbolic surfaces} or, more generally, on closed \emph{hyperbolic 2-orbifolds}. 
We refer to these two types as \term{classical Anosov flows}. 
They belong to the larger family of \term{algebraic Anosov flows}, which consists of those Anosov flows that can be obtained as the quotient by a cocompact subgroup of a Lie group where a 1-parameter family is acting, and were characterized by Tomter~\cite{Tomter}. 
In addition to the two previous families, there exist entire families of Anosov flows that are not algebraic; these are constructed using a collection of \emph{3-manifold cut-and-paste} techniques, usually referred to as \emph{surgeries} \cites{FW, HandelThurston, Goodman_surgery} or \emph{Lego-type constructions} \cites{BBY, Paulet}.

Two flows are \term{orbitally equivalent} if there exists a homeomorphism between the underlying manifolds that maps orbits of the first flow onto orbits of the second, preserving the natural orientation of the orbits induced by the flow action. 
The question of whether two Anosov flows are orbitally equivalent can be answered for some special families of Anosov flows. For instance, for the family of classical Anosov flows: two geodesic flows are equivalent if and only if the underlying 2-orbifolds are of the same type (cf. \cite{Gromov}); two suspensions are equivalent if and only if the underlying matrices are conjugate in~$\GLZ$ (cf. \cite{Plante}); and a geodesic flow is never equivalent to a suspension. However, the general problem of classifying Anosov flows up to orbital equivalence remains a fundamental open problem in the study of Anosov flows. 

Following the constructions of Goodman \cite{Goodman_surgery} and Fried \cite{FriedAnosov} elaborating on the notion of Dehn surgery, a more flexible notion of equivalence was proposed by Christy and Ghys~\footnote{This notion appears in Kirby’s list~\cite{Kirby} (Problem 3.54, proposed by Christy), a reference brought to our attention by Bin Yu and Chi Cheuk Tsang. It traces back to a oral questions by Ghys in the 90s.}: 
two flows are \term{almost equivalent} if there exists a homeomorphism from the complement of a finite number of periodic orbits of the first flow to the complement of the same number of periodic orbits of the second flow, which maps orbits onto orbits while preserving their orientation. A seminal construction by Birkhoff \cite{Birkhoff} and Fried \cite{FriedAnosov} shows that the geodesic flow on a hyperbolic surface is almost equivalent to the suspension flow of an automorphism of the torus, and that automorphism has been calculated explicitly~\cites{GhysGV, Hashiguchi_first-return}.

Christy and Ghys proposed the general problem of understanding when two Anosov flows are almost equivalent. It should be noted that Anosov flows can be \emph{transitive} (i.e., there exists a dense trajectory in the supporting 3-manifold) or not, and their invariant \emph{center-stable} and \emph{center-unstable} foliations can be \emph{orientable} or not. Examples exist for all four possibilities. Since almost equivalence between Anosov flows preserves each of these properties, there are many distinct classes under this relation. A relevant problem is thus the following:

\begin{problem}[Christy-Ghys]
\label{christy-ghys_question}
Determine if all the Anosov flows defined on closed orientable 3-manifolds that are: 
\begin{enumerate}[(i)]
    \item transitive,
    \item have orientable center-stable and center-unstable foliations,
\end{enumerate}
are pairwise almost equivalent.     
\end{problem}

Classical Anosov flows always satisfies (i) above, and the subset satisfying (ii) are precisely the suspensions of hyperbolic automorphisms of the 2-torus with positive trace, and geodesic flows on hyperbolic orientable orbifolds. 
Problem \ref{christy-ghys_question} was answered positively by Minakawa for the family of suspension Anosov flows. 

\begin{introtheorem}[Minakawa~\cite{Minakawa}]\label{T:Main}
Every Anosov flow that is the suspension of a hyperbolic automorphism of the 2-torus with positive trace is almost equivalent to the suspension flow generated by the matrix~$(\begin{smallmatrix}2&1\\1&1\end{smallmatrix})\in\SLZ$.
\end{introtheorem}

A weak version of this statement had already been proven~\cite{Commensurability}, where almost equivalence was replaced by \emph{almost commensurability} (allowing for finite coverings). The family treated by Minakawa corresponds to the suspension of matrices $A\in\SLZ$ with $\tr(A)\geq 3$. 
For such a matrix, we denote by
\begin{itemize}
    \item $\TT^3_A$ the 3-manifold~$\TT^2\times[0,1]/_{(x,1)\sim(Ax,0)}$

    \item$\Phisus$ the flow on~$\TT^3_A$ tangent to the $[0,1]$-coordinate, called the \emph{suspension flow}.  
\end{itemize}
Theorem~\ref{T:Main} is implied by the following more precise result:

\begin{introtheorem}[Minakawa~\cite{Minakawa}]
\label{T:Trace}
If $A\in\SLZ$ has a trace greater than 3, then there exists~$B\in\SLZ$ with $3\le \tr B < \tr A$ such that the suspension flow on~$\TT^3_A$ is almost equivalent to the suspension flow on~$\TT^3_B$.
\end{introtheorem}

The first purpose of this note is to provide a written proof of Minakawa's Theorems ~\ref{T:Main} and~\ref{T:Trace}, which is a very important step regarding Problem \ref{christy-ghys_question}, and for which there is no available written reference\footnote{The reference ~\cite{Minakawa} is a video of a talk given in Tokyo in 2013 where Minakawa announced the theorem. The full proof is not in the video; it was outlined in the abstract and is currently not available online. H. Minakawa informed us in a personal communication that the result was already announced in 2004 at the 51st Topology Symposium at Yamagata.}. The second purpose is to improve Minakawa's result a bit further, as explained below.

Define the \term{graph of Anosov flows}, or \term{Ghys graph},  to be the graph $\ggood$ such that 
\begin{enumerate}
    \item The vertices are the orbital equivalence classes of Anosov flows on 3-manifolds;

    \item Two vertices are connected by an edge if one can remove one periodic orbit from each flow and obtain two orbitally equivalent flows.
\end{enumerate}
As noted above, this graph has several connected components (corresponding to the transitivity and orientability properties of the invariant foliations). By treating this graph as a metric space where edges have length 1, the \term{distance} $\dgood$ between two Anosov flows is the minimal number of periodic orbits to be removed from both flows to obtain an almost equivalence (or infinite if the flows are not almost equivalent). 

Theorem~\ref{T:Main} can be rephrased as follows: \emph{All suspension Anosov flows on orientable 3- manifolds with orientable center-stable/unstable foliations lie in the same connected component of~$\ggood$}. 
Our improvement of Theorems \ref{T:Main} and \ref{T:Trace} consists in estimating distances in the graph $\ggood$. 

Let $X$ denote the matrix~$(\begin{smallmatrix}1&1\\0&1\end{smallmatrix})$ and $Y$ the matrix~$(\begin{smallmatrix}1&0\\1&1\end{smallmatrix})$. If we choose a hyperbolic matrix $W$ with positive entries, then the matrix $ZW$ is also hyperbolic for either choice of $Z=X$ or $Z=Y$. This follows simply by checking that the trace of $ZW$ is greater than $2$:
\begin{align*}
    & \begin{pmatrix}
        1 & 1\\
        0 & 1
    \end{pmatrix}
    \begin{pmatrix}
        w_{11} & w_{12}\\
        w_{21} & w_{22}
    \end{pmatrix}=
    \begin{pmatrix}
        w_{11}+w_{21} & w_{12}+w_{22}\\
        w_{21} & w_{22}
    \end{pmatrix}\Longrightarrow \tr(XW)=\tr(W)+w_{21}\geq\tr(W),\\
    & \begin{pmatrix}
        1 & 0\\
        1 & 1
    \end{pmatrix}
    \begin{pmatrix}
        w_{11} & w_{12}\\
        w_{21} & w_{22}
    \end{pmatrix}=
    \begin{pmatrix}
        w_{11} & w_{12}\\
        w_{21}+w_{11} & w_{22}+w_{12}
    \end{pmatrix}\Longrightarrow \tr(YW)=\tr(W)+w_{12}\geq\tr(W).
\end{align*}
Since the entries of $W$ are positive and the matrix is hyperbolic then necessarily $\tr(W)\geq 3$, which implies that $ZW$ is also hyperbolic. The next result states that the vertices in the graph $\ggood$ corresponding to the orbital equivalence classes of $(\TT^3_{W},\Phisus)$ and $(\TT^3_{ZW},\Phisus)$ are in the same connected component, and a bound on their distance in the metric $\dgood$ is provided.

\begin{introtheorem}\label{T:Distance}
Let $W$ be a matrix in $\SLZ$ with positive entries and $\tr(W)\geq 3$. Then it is verified  
\[\dgood((\TT^3_{W},\Phisus), (\TT^3_{ZW}, \Phisus))\le 3, \quad \text{where } Z=X \text{ or } Z=Y.\]
\end{introtheorem}

Let us see how to derive some metric properties of the graph $(\ggood,\dgood)$ using the previous theorem. We say that a 3-manifold is a \emph{sol-manifold} if its fundamental group is solvable. Let $\ggood_\circ$ denote the set of vertices of $\ggood$ representing orbital equivalence classes of Anosov flows defined on sol-manifolds, and whose center-stable and center-unstable foliations are orientable. Using Theorem \ref{T:Distance} we can give bounds on the distance $\dgood$ between elements in $\ggood_\circ$.  

In \cite{Plante} it is shown that every Anosov flow on a sol-manifold is orbitally equivalent to a suspension Anosov flow $(\TT^3_{A},\Phisus)$, where $A$ is some hyperbolic matrix in $\SLZ$. Moreover, recall that two suspension Anosov flows $(\TT^3_{A},\Phisus)$ and $(\TT^3_{B},\Phisus)$ are orbitally equivalent if and only if the matrices $A$ and $B$ are conjugated in $\GLZ$. Therefore, the elements in $\ggood_\circ$ can be indexed over the set of $\GLZ$-conjugacy classes of hyperbolic matrices in $\SLZ$ with positive trace. We will use the following:

\paragraph{\textbf{Representing hyperbolic matrices with words in two symbols}}- 
It is a \emph{folklore result} that every hyperbolic matrix with positive trace is conjugated to a positive product of the matrices $X$ and $Y$ containing both letters, and this product is unique up to cyclic permutation (see, e.g., Prop. 4.3 of~\cite{Lorenz}). Therefore, the elements of $\ggood_\circ$ can be indexed over the set of finite positive words on the alphabet $\{X,Y\}$ up to cyclic permutation. Observe that, up to conjugation, there is only one hyperbolic matrix with trace equal to 3 (the minimum possible value of the trace) and this is the matrix $(\begin{smallmatrix}2&1\\1&1\end{smallmatrix})=XY$. 

Following the previous discussion, any element of $\ggood_\circ$ can be represented by a suspension Anosov flow $(\TT^3_{W},\Phisus)$, where the matrix $W$ equals a (unique up to cyclic permutation) product $W=Z_1^{m_1}\cdots Z_k^{m_k}$ such that each $Z_i$ equals either $X$ or $Y$, every $m_i$ is positive, and there is at least one $X$ and one $Y$ in the decomposition. Assume that the \emph{word length} $m=\sum_{i=1}^k m_i$ of $W$ is much bigger that $2$, and write $W=Z_1W_1$, where either $Z_1=X$ or $Z_1=Y$ and $W_1=Z_1^{m_1-1}\cdots Z_k^{m_k}$ is a (smaller) positive word in $\{X,Y\}$ containing both letters. Since $W_1$ is also hyperbolic and all its entries are positive, then Theorem \ref{T:Distance} says that the class of the suspension flow $(\TT^3_{W},\Phisus)$ is at distance at most 3 from the class of $(\TT^3_{W_1},\Phisus)$. Hence, the effect of deleting a letter on the left of the word $W$ produces a new suspension Anosov flow (in the class of the smaller word $W_1$) at distance at most $3$ from the suspension of $W$. By deleting letters either on the left or on the right of the word $W$, we can construct a explicit paths of almost equivalences starting at the class of $(\TT^3_W,\Phisus)$ and ending at the class of $(\TT^3_{XY},\Phisus)$, and we see that 
$$\dgood((\TT^3_{XY},\Phisus),(\TT^3_{W},\Phisus))\leq 3m-6.$$

The previous bound can be improved. Given a hyperbolic matrix $A$ in $\SLZ$ with positive trace, define the \emph{word length} $m(A)$ of $A$ to be the word length in the alphabet $\{X,Y\}$ of the matrix $W$, the (unique up to cyclic permutations) positive product of $X$ and $Y$ that is conjugated to $A$. We will prove\footnote{We thank the referee for suggesting the bound involving matrix traces.}:

\begin{corollary}\label{C:Distance}
For every matrix $A$ in $\SLZ$ with $\tr(A)\geq 3$, the distance in $\ggood$ from the suspension of $A$ to the suspension of $XY=(\begin{smallmatrix}2&1\\1&1\end{smallmatrix})$ verifies the following bounds:
\begin{enumerate}[(i)]
    \item $\dgood((\TT^3_{XY},\Phisus),(\TT^3_{A},\Phisus))\leq 2m(A)-4$,

    \item $\dgood((\TT^3_{XY},\Phisus),(\TT^3_{A},\Phisus))\leq 2\tr(A)-5.$ 
\end{enumerate}
\end{corollary}

\medskip

Theorem~\ref{T:Distance} can be rephrased in terms of \emph{Birkhoff sections}, which is a notion that provides the adequate setting to work with almost equivalences. Recall that a \emph{global transverse section} for a non-singular flow is an embedded closed surface in the supporting manifold that transversely intersects every orbit segment within a uniformly bounded time. A Birkhoff section is a generalization of this notion, allowing the surface to have finitely many boundary components where it fails to be transverse to the flow. 

\begin{definition}\label{S:BirkhoffSection}
Given a flow~$\Phi=(\Phi^t)_{t\in\RR}$ on a compact 3-manifold~$M$, a~\term{Birkhoff surface} is an immersion~$\iota:(S,\partial S)\to (M,\Gamma)$ of a compact surface $S$ (possibly with boundary) inside $M$, taking the boundary components $\partial S$ onto a finite set $\Gamma$ of periodic orbits of $\Phi$, such that: 
\begin{enumerate}[(i)]
\item restricted to the interior $\textrm{int}(S):=S\setminus\partial S$ the map $\iota$ is an embedding and the image $\iota(S\setminus\partial S)$ is positively transverse to the $\Phi$-orbits;
\item restricted to any component of $\partial S$ the map $\iota$ is a covering onto a periodic orbit in $\Gamma$.
\end{enumerate}
The surface~$S$ is a \term{Birkhoff section} if it further satisfies the condition:
\begin{enumerate}
\item[(iii)] there exists $T>0$ such that every $\Phi$-orbit segment of length $T$ intersects $\iota(S)$. 
\end{enumerate}
\end{definition}

When discussing about Birkhoff sections, we will in general omit talking about the immersion $\iota$ and we will identify $S$ directly with its image $\iota(S)$ inside $M$. Observe that a global transverse section can be seen as the special case of a Birkhoff section having empty boundary. Properties (i), (ii), (iii) above imply the existence of an induced \term{first-return map} $f_S:\textrm{int}(S)\to\textrm{int}(S)$, only defined on the interior of $S$. Consequently, the flow $\Phi$ restricted to the open 3-manifold $M\setminus\Gamma$ has $\textrm{int}(S)$ as global transverse section, and is equivalent to the suspension flow induced by $f_S$. 

Associated to a Birkhoff section there is an induced homeomorphism $\hat{f}_S$, defined on the closed surface $\widehat{S}$ obtained by collapsing each boundary component of $S$ into a point. These collapsed components become periodic orbits of $\hat{f}_S$. 
This homeomorphism is called the \term{blow-down} associated to the Birkhoff section. 
In the case where $\Phi$ is an Anosov flow with orientable center-stable and center-unstable foliations it can be shown that is a {pseudo-Anosov} homeomorphism with even prongs only (cf. \cite{FriedAnosov}).

The key observation is that, if $\Phi$ has a Birkhoff section $\iota:(S,\partial S)\to (M,\Gamma)$, then by removing from $M$ the periodic orbits of~$\Phi$ that form the image $\Gamma=\iota(\partial S)$, we obtain an almost equivalence between~$\Phi$ and the suspension flow induced by~$\hat{f}_S:\widehat{S}\to\widehat{S}$. Therefore, Theorem~\ref{T:Distance} is a consequence of the following:

\begin{introtheorem}
\label{T:DistanceBis}
Let $W$ be a matrix in $\SLZ$ with positive entries and $\tr(W)\geq 3$. Then the Anosov flow $(\TT^3_{ZW},\Phisus)$, where $Z=X$ or $Z=Y$, admits a genus-one Birkhoff section with at most three boundary components, and whose induced first-return map is given by the matrix~$W$.
\end{introtheorem}

The topological or metrical properties of the graph $\ggood$ have been mostly unexplored until the writing of this note. 
We end this introduction with a small survey of known results related to Problem \ref{christy-ghys_question} and the topology of this graph, as well as some related questions. 

\paragraph{\textbf{Minimal-genus Birkhoff sections}}
Fried proved \cite{FriedAnosov} that every transitive Anosov flow admits (infinitely many) Birkhoff sections. Following Minakawa's theorem above, Problem \ref{christy-ghys_question} reduces to the following problem, originally stated by Fried~\cite{FriedAnosov}:
\begin{problem}[Minimal genus problem]
\label{min-genus_question}
Determine if all the Anosov flows defined on closed orientable 3-manifolds that are: 
\begin{enumerate}[(i)]
    \item transitive,
    \item have orientable center-stable and center-unstable foliations,
\end{enumerate}
admit a genus-one Birkhoff section.     
\end{problem}
\noindent
Because, if two such flows are endowed with genus-one Birkhoff sections, then each one is almost equivalent to a suspension and, since by Theorem \ref{T:Main} all the suspensions are pairwise almost equivalent, if follows that the original flows are almost equivalent. 

\paragraph{\textbf{Geodesic flows}}
Several works are devoted to describe genus-one Birkhoff sections for geodesic flows on hyperbolic surfaces and hyperbolic 2-orbifolds; see for instance \cites{FriedAnosov, GhysGV, Hashiguchi_first-return, Brunella, HM, GenusOne}. Together with Minakawa's Theorem~\ref{T:Main}, these results provide several examples of Anosov geodesic flows that are almost equivalent to the suspension flow generated by the matrix~$(\begin{smallmatrix}2&1\\1&1\end{smallmatrix})$. Moreover, some of the constructions in the articles cited above give an explicit calculation of the first return map onto the genus-one Birkhoff section, allowing to obtain bounds for the distance $\dgood$ between some geodesic flows and and the suspension of~$(\begin{smallmatrix}2&1\\1&1\end{smallmatrix})$. 
However, at the moment of writing this note, \emph{we do not know whether all Anosov geodesic flows admit genus-one Birkhoff sections or not}.
\footnote{
In a previous version of this note we stated that all classical Anosov flows admit a genus-one Birkhoff section. However, we found a gap in our argument during the revision process. The cases of geodesic flows where we do have explicit genus-one Birkhoff sections are when the orbifold has genus 0, or when it has positive genus $g$ and at most $2g+6$ conic points. 
}
\begin{question}
    Does every geodesic flow on a hyperbolic 2-orbifold admit a genus-one Birkhoff section?
\end{question}

More generally one can wonder about the minimal genus for Anosov flows on Seifert fibered manifold. 
By Ghys' theorem~\cite{GhysAnosov} these flows are all finite coverings of geodesic flows on hyperbolic 2-orbifolds, and it was recently shown~\cite{BF} that for every admissible degree of such a covering, there is only one or two Anosov flows, up to orbital equivalence. 

\begin{question}
 Does every finite covering of a geodesic flow admit a genus-one Birkhoff section? 
\end{question}

\paragraph{\textbf{Graph manifolds}}
In \cite{Tsang2} it is shown that every Anosov flow on a \emph{graph manifold} $M$ that is \emph{totally periodic} and satisfies the conditions of Problem \ref{min-genus_question} admits a genus-one Birkhoff section, and hence is almost equivalent to a suspension.

\paragraph{\textbf{Suspensions of Penner maps}}
In \cite{Tsang2} it is also shown that every Anosov flow that admits a Birkhoff section (of arbitrary genus) whose first-return map is of Penner type with orientable stable and unstable foliations admits a genus-one Birkhoff section, and hence is almost equivalent to a suspension.

\paragraph{\textbf{Birkhoff sections with minimal number of boundary components}} 
In \cite{Marty} it is shown that every Anosov flow that is $\RR$-covered and satisfies the hypothesis of Problem \ref{min-genus_question} admits a Birkhoff section with only one boundary component (although there is no control on the genus of the surface), and in~\cite{Tsang1} it is shown that for transitive non-$\RR$-covered flows this number is two. 

\paragraph{\textbf{Almost equivalence and orbit space}}
In \cite{Bonatti-Iako} (and other related works as \cite{Fenley_Anosov-flows}) it is studied how the orbit space of Anosov flows changes under almost equivalences or, equivalently, under \emph{Fried surgeries} (see definitions in the referred work). For example, it is shown that every Anosov flow in the hypothesis of Problem \ref{christy-ghys_question} is almost equivalent to a $\RR$-covered Anosov flow. 

\paragraph{\textbf{Distances in the graph} $\ggood$} 
Let $\ggood_\textrm{sus}$ denote the subset of vertices of $\ggood$ that correspond to suspensions of hyperbolic matrices with positive trace. Theorem~\ref{T:Distance} and Corollary~\ref{C:Distance} only give a coarse bound for the distance between an arbitrary element of $\ggood_\textrm{sus}$ and $(\TT^3_{XY},\Phisus)$ (the simplest suspension Anosov flow). 

\begin{question}
Is it possible to give a more accurate metric description of $(\ggood_\textrm{sus},\dgood)$ ? 
\end{question}

For instance, consider the graph $\mathrm{W}$ whose vertices are classes of positive words on $\{X,Y\}$ identified up to cyclic permutations, and where two vertices are connected with one edge if we can go from a word representing a vertex (up to cyclic permutation) to the other by adding a letter to the left. \emph{Which are the metric properties of the embedding $(\mathbf{W},\mathrm{dist}_\mathbf{W})\to(\ggood_\mathrm{sus},\dgood)$} ?

From a more general point of view, it would be desirable to have a description of the large-scale topology of the graph of Anosov flows. The graph $\ggood$ allows to represent surgery operations between Anosov flows (see, e.g. \cite{Sha24}), and a large-scale metric description $\ggood$ would reflect the complexity in the surgery process to obtain a given Anosov flow from a simpler one.   


\paragraph{\textbf{Loops in the graph} $\ggood$}
In \cite{Sha24} it is shown that, at some specific vertices of $\ggood$ corresponding to suspensions, there exist infinitely many different closed loops of length 2. Moreover, an asymptotic estimation on the number of these loops is given. This is done by showing that every suspension Anosov flow admits infinitely many genus-one Birkhoff sections, having exactly two boundary components. 
On the other hand, we do not know whether there exist \emph{Anosov cosmetic surgeries}, {\it i.e.}, loops of length~1 in~$\ggood$. 

\bigskip
\noindent{\bf Organization of the paper.}
In Section~\ref{S:Fried}, we describe an operation on transverse surfaces called \emph{Fried sum}, that is needed along the main construction. In Section~\ref{S:Minakawa}, we explain Minakawa's construction of a genus-one Birkhoff section for a suspension Anosov flow, and prove Theorem~\ref{T:DistanceBis}.

\section{Fried sum and Euler characteristic}
\label{S:Fried}

In this section we describe an operation introduced by Fried~\cite{FriedAnosov} that takes two Birkhoff surfaces and produces a new one. 

\subsection{Fried sum of two Birkhoff surfaces}\label{S:FriedSum}
Let $M$ be a closed 3-manifold, $\Phi$ a flow generated by a non-singular vector field $X$ on $M$, and $S_1, S_2$ two Birkhoff surfaces. Their boundaries are links $\Gamma_1, \Gamma_2$ in $M$ (not necessarily disjoint) consisting of unions of periodic orbits of $\Phi$. We denote the union of these links by $\Gamma = \Gamma_1 \cup \Gamma_2$. 

By perturbing the surfaces while keeping transversality with $X$ along their interiors and keeping $\Gamma$ invariant, one may assume that $S_1$ and $S_2$ are in transverse position. Their intersection $S_1 \cap S_2$ is thus a 1-manifold, consisting in a union of circles and arcs with endpoints in~$\Gamma$.  
From the union $S_1 \cup S_2$, which is not an embedded surface when $S_1\cap S_2$ is non-empty, we can construct a new Birkhoff surface by modifying the union in a neighbourhood of $S_1 \cap S_2$. This process is known as \emph{Fried's desingularization}.

\begin{figure}[t]
\begin{center}
\includegraphics[width=0.8\textwidth]{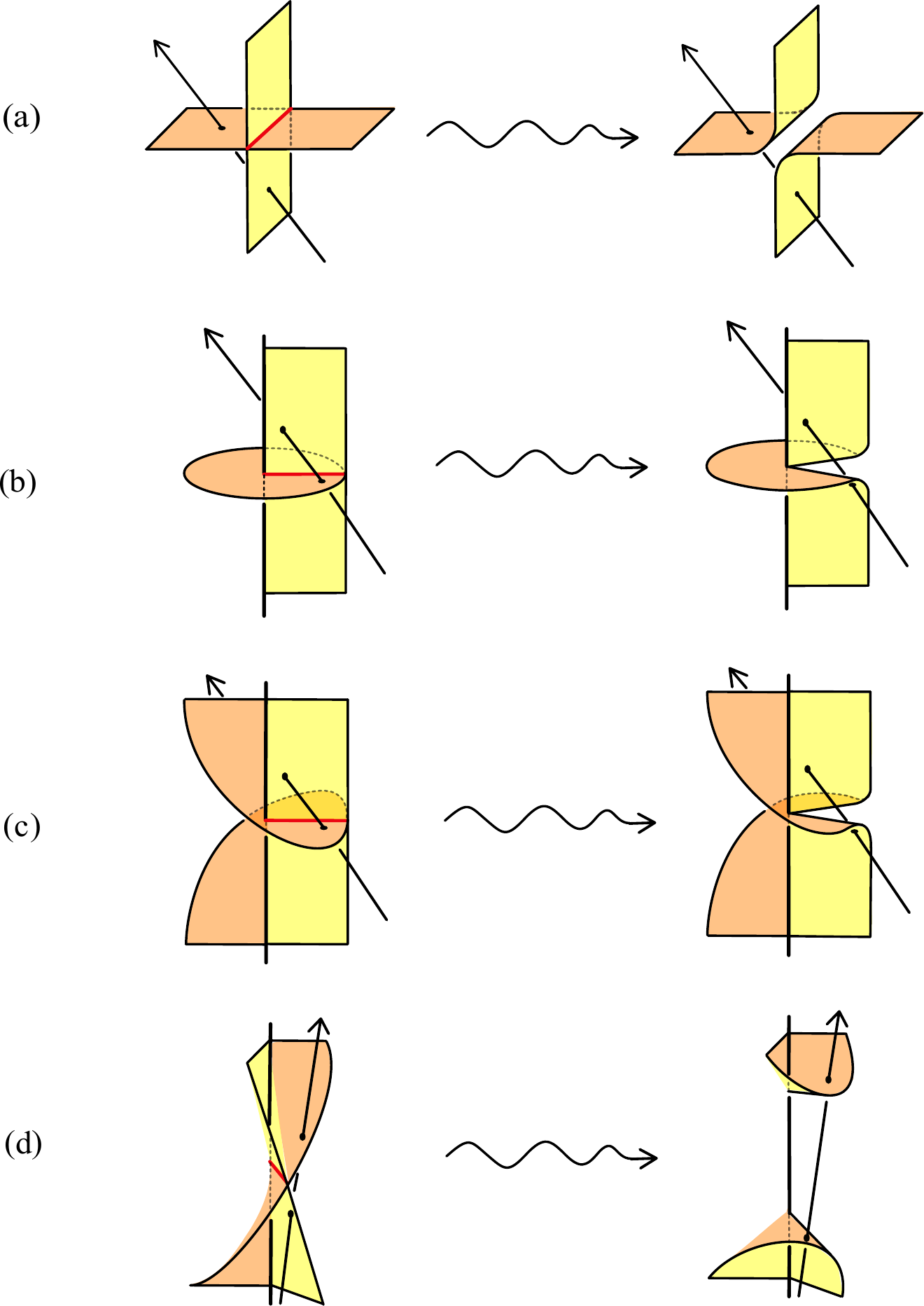}
\caption{Fried desingularization.}
\label{F:DesArc}
\end{center}
\end{figure}

Let $X$ be the closure of the set $S_1\cap S_2\setminus\Gamma$, which is a 1-manifold consisting of circles and arcs with endpoints in $\Gamma$. Observe that, in the case where $\Gamma_1\cap\Gamma_2=\emptyset$ then $X=S_1\cap S_2$. When $\Gamma_1$ and $\Gamma_2$ share a component $\gamma$, the set $X$ excludes from $S_1 \cap S_2$ all the points in $\gamma$ that are not accumulated by interior points of both $S_1$ and $S_2$. 
Let $x$ be a point in $X$ and $U$ a small open neighbourhood around $x$. The desingularization at $x$ consists in removing $(S_1\cup S_2)\cap U$ from the union $S_1\cup S_2$ and gluing back one of the following local model surfaces; in each case, there are two possibilities for making the desingularization of the surface; however, only one of them (up to isotopy) produces a surface transverse to $X$:

\begin{itemize}
    \item If $x \notin \Gamma$, we use the model shown in Figure\eqref{F:DesArc}-(a);
    \item If $x \in \Gamma$ and $x$ is in the interior of $S_1$ or $S_2$, we use the model shown in Figure\eqref{F:DesArc}-(b);
    \item If $x \in \Gamma$, $x\in \partial S_1$, $x\in \partial S_2$, and the orientations locally agree, we use the model shown in Figure\eqref{F:DesArc}-(c);
    \item If $x \in \Gamma$, $x\in \partial S_1$, $x\in \partial S_2$, and the orientations locally disagree, we use the model shown in Figure\eqref{F:DesArc}-(d); in this case, the operation is only semi-local: in order for the new boundary component of the surface to be a covering of the considered orbit, one has to isotope the surface in a tubular neighbourhood of the orbit.
\end{itemize}

Cases (c) and (d) do not show up in this note, but we include them for completeness. 

\begin{definition}\label{D:FriedSumDef}
Given two Birkhoff surfaces $S_1, S_2$ as above, their \term{Fried sum}, denoted by $S_1 \scup S_2$, is the surface obtained from $S_1 \cup S_2$ by desingularizing all circles and arcs of $S_1 \cap S_2$ transversely to the vector field $X$.
\end{definition}

\begin{remark}
It is possible to check that the result $S_1\#S_2$ of desingularizing two Birkhoff surfaces $S_1$ and $S_2$ is another Birkkhoff surface, in the sense that $S_1\scup S_2$ is the image of an immersion $\Sigma\to M$ that is an embedding transverse to the flow along $\textrm{int}(\Sigma)$, and a covering onto a finite set of periodic orbits along $\partial\Sigma$. It is also possible to check that the immersion $\Sigma\to M$ obtained by desingularizing $S_1\#S_2$ is unique up to isotopies obtained by pushing along flow-lines. This last fact can be checked using the techniques at Section~$3$ of \cite{Shannon}. If $\iota_i:\Sigma\to M$, $i=1,2$ are two Birkhoff sections obtained by desingularizing $S_1\# S_2$, then both surfaces have the same boundary link $\Gamma$ and for each component $\gamma\in\Gamma$ they have the same combinatorial invariants (cf. \emph{linking number} and \emph{multiplicity} at \cite{Shannon}). Applying Proposition 3.17 from \cite{Shannon} at each component of $\Gamma$, we can see that both embeddings are isotopic pushing along flow-lines.  
\end{remark}

\subsection{Normal blow-up along $\Gamma$}\label{S:NormalBlowUp}
The Fried desingularization of $S_1$ and $S_2$ can be better understood by blowing up the 3-manifold $M$ along the link $\Gamma:=\Gamma_1\cup\Gamma_2$, formed by the union of the boundary links $\Gamma_i$ of each surface $S_i$. 

The \emph{normal blow-up} of $M$ along $\Gamma$ is the space $M_\Gamma$ obtained by replacing each point $p \in \Gamma$ by its normal sphere bundle $S((TM)_p / \mathbb{R} X(p))$, which is topologically a circle. With the adequate topology, $M_\Gamma$ is a 3-manifold whose boundary $\partial M_\Gamma$ is a union of tori (one for each component of $\Gamma$) and is a compactification of $M \setminus \Gamma$. The flow $\Phi$ extends to a flow $\Phi^*$ on $M_\Gamma$ that preserves its boundary. The surfaces $S_1$ and $S_2$ extend to embedded surfaces $(S_1^*, \partial S_1^*)$ and $(S_2^*, \partial S_2^*)$ in $(M_\Gamma, \partial M_\Gamma)$ that are transverse to the $\Phi^*$-orbits and transverse to each other. After perturbation, we may assume these conditions also hold on their boundaries within $\partial M_\Gamma$. 

\begin{lemma}\label{lemma_desingularizing_boundary}
The surface $S_1^*\#S_2^*$ obtained by desingularizing the union of $S_1^*$ and $S_2^*$ in $M_\Gamma$, is homeomorphic to the image of $S_1 \scup S_2$ under normal blow-up along~$\Gamma$. 
\end{lemma}

\begin{proof}
The proof follows directly from the fact that both surfaces are obtained by desingularizing homeomorphic sets in $M_\Gamma \setminus \partial M_\Gamma$ and $M \setminus \Gamma$. 
\end{proof}

Lemma \ref{lemma_desingularizing_boundary} is particularly useful for understanding the embedding $\Sigma\to(S_1 \scup S_2)$, since $\Sigma$ is homeomorphic to the blown-up surface $(S_1\scup S_2)^*$, and it is visually simpler to perform the desingularization of the boundaries $\partial S_1^* \cup \partial S_2^*$ within each component of $\partial M_\Gamma$. 

\subsection{Euler characteristic}\label{S:Euler}
Another application of Lemma \ref{lemma_desingularizing_boundary} is to determine the Euler characteristic of $\Sigma$ in the embedding $\Sigma\to S_1 \scup S_2$, which is needed during the proof of Theorem \ref{T:DistanceBis} in order to check that the obtained Birkhoff section has genus 1.

\begin{lemma}\label{lemma_linear_Euler_char}
With the notations above, $\chi(S_1 \scup S_2):=\chi(\Sigma) = \chi(S_1^*) + \chi(S_2^*)$. 
\end{lemma}

\begin{proof}
Note that to construct the surfaces $S_i^*$ one must remove a disk from $S_i$ for each intersection of $\partial S_j$ ($j \neq i$) with the interior of $S_i$, to ensure that $S_1^*$ and $S_2^*$ are properly defined in $M_\Gamma$. To prove the lemma, triangulate both surfaces $S_i^*$ in such a way the intersection $S_1^*\cap S_2^*$ lies in the 1-skeleton of both triangulations. In the blown-up manifold, the desingularization process effectively "uncouples" the surfaces along $S_1^*\cap S_2^*$ by doubling all 0- and 1-simplices contained in this intersection. Therefore, the Fried sum $S_1^*\scup S_2^*$ in $M_\Gamma$ can be triangulated with $n_k$ simplices for each dimension $k=0,1,2$, where $n_k$ equals the sum of the $k$-simplices in $S_1^*$ and $S_2^*$. It follows that $\chi(S_1^* \scup S_2^*) = \chi(S_1^*) + \chi(S_2^*)$. Since $\Sigma$ is homeomorphic to the blown-up surface $(S_1\scup S_2)^*$, by Lemma \ref{lemma_desingularizing_boundary} we conclude that $\chi(\Sigma)=\chi(S_1^*)+\chi(S_2)^*$. 
\end{proof}

\section{Proof of the main results}\label{S:Minakawa}
Assume we are given a hyperbolic matrix $W$ in $\SLZ$ with positive entries. Let $Z$ be either $Z=X$ or $Z=Y$, where the matrices $X, Y$ are given by  
$$X=
\begin{pmatrix}
    1 & 1\\
    0 & 1
\end{pmatrix},\ 
Y=
\begin{pmatrix}
    1 & 0\\
    1 & 1
\end{pmatrix}.
$$
Since $W$ is hyperbolic and has positive entries then $ZW$ is also hyperbolic with positive entries. Consider the manifold $\TT^3_{ZW}$ equipped with its corresponding suspension flow $\Phisus$. Our goal is to construct a genus-one Birkhoff section for this flow whose first-return map is conjugated to the matrix~$W$. 

This Birkhoff section will be obtained by performing a Fried sum (see Section~\ref{S:Fried}) between a global transverse section~$\Th := \TT^2 \times \{2/3\}$, and a Birkhoff surface $\PP$ that is an immersion of a pair of pants and whose boundary link consists of three periodic orbits. 

There are two main points to verify: 
\begin{enumerate}
    \item The surface obtained via the Fried sum has genus 1. 
    \item The first-return map is indeed given by the matrix~$W$ in an appropriately chosen basis. 
\end{enumerate}

Item (1) above requires a careful choice of pair of pants~$\PP$, since most configurations would result in higher-genus sections. The construction of this Birkhoff section and the calculation of its first-return map are detailed in Sections~\ref{S:The_pair_of_pants}, \ref{S:Cutting}, \ref{S:FirstReturn}, and \ref{S:case_YW}. The proofs of Theorems~\ref{T:Trace}, \ref{T:Distance}, and \ref{T:DistanceBis} are then completed in Section~\ref{S:proof_thmB}.
Since the proofs for both cases $Z=X$ and $Z=Y$ are analogous, we will focus on the case $Z=X$ in the following sections and describe the necessary modifications for the case $Z=Y$ in Section~\ref{S:case_YW}. 


\subsection{The pair of pants}
\label{S:The_pair_of_pants}
Let $W$ be a hyperbolic matrix with positive entries. Consider the matrix $XW$ and denote its entries by $\left(\bsm a&b\\c&d\esm\right)$. Define $t := \tr(XW) = a+d$. 
We describe first the construction of a a parallelogram contained in~$\TT^2$, whose vertices are fixed points of the action of~$XW$, and whose sides consist of specifically chosen segments. This parallelogram will be used later to construct an embedded pair of pants in the suspension manifold~$\TT^3_{XW}$.

\begin{figure}[t]
\centering
\includegraphics[width=0.9\textwidth]{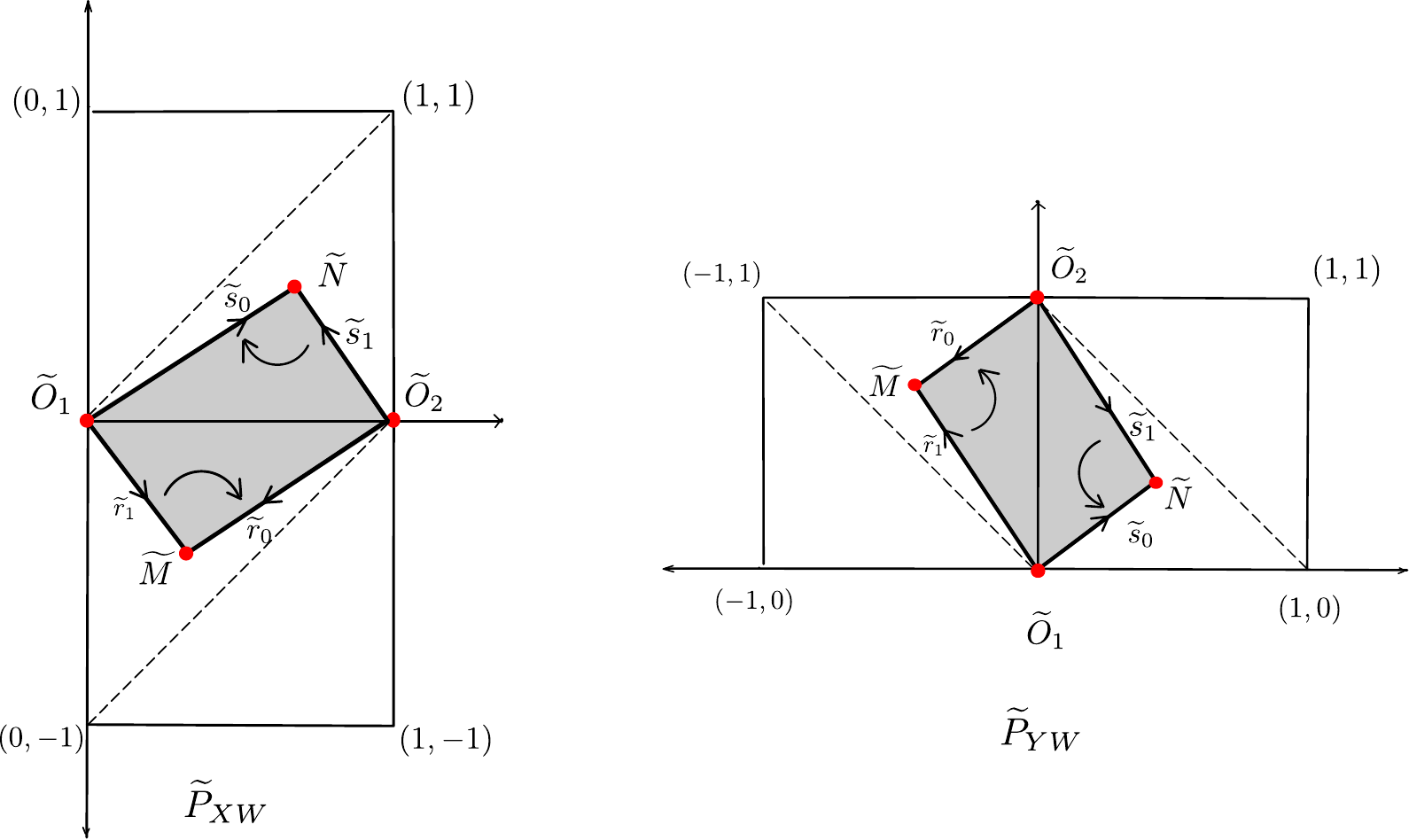}
\caption{The polygons $\widetilde{P}_{XW}$ (left) and $\widetilde{P}_{YW}$ in $\RR^2$.}
\label{pantalon-r2}
\end{figure}

\paragraph{\textbf{Construction of the parallelogram}}
Consider the vector in $\RR^2$ given by 
$$v = 
v\begin{pmatrix}
    a & b\\
    c & d
\end{pmatrix} := 
\begin{pmatrix} 
(d-1)/(t-2) \\ -c/(t-2) 
\end{pmatrix}.$$ 
A direct calculation shows that $XW \cdot v = v - (1,0)$. That is:
\begin{equation}\label{eq_minakawa_vector}
\bpm a&b\\c&d\epm \bpm (d-1)/(t-2) \\ -c/(t-2) \epm = \bpm (1-a)/(t-2) \\ -c/(t-2) \epm = \bpm (d-1)/(t-2) \\ -c/(t-2) \epm - \bpm 1\\0\epm.
\end{equation}
Using the vector $v$, we can construct a parallelogram $\widetilde{P}=\widetilde{P}_{XW}$ in $\RR^2$ with the following sides and vertices (see left-side of Figure~\eqref{pantalon-r2}):
\begin{align*}
&\tilde{r}_1 = \{t \cdot v : 0 \le t \le 1\}, &&\widetilde{O}_1 = (0,0),\\
&\tilde{r}_0 = XW \cdot \tilde{r}_1 + (1,0), &&\widetilde{M} = \left(\frac{d-1}{t-2}, \frac{-c}{t-2}\right),\\
&\tilde{s}_1 = -\tilde{r}_1 + (1,0), &&\widetilde{O}_2 = (1,0),\\
&\tilde{s}_0 = -XW \cdot \tilde{r}_1, &&\widetilde{N} = \left(\frac{a-1}{t-2}, \frac{c}{t-2}\right).  
\end{align*}

\begin{figure}[t]
\centering
\begin{picture}(160,75)(0,0)
\put(0,0){\includegraphics[width=.42\textwidth]{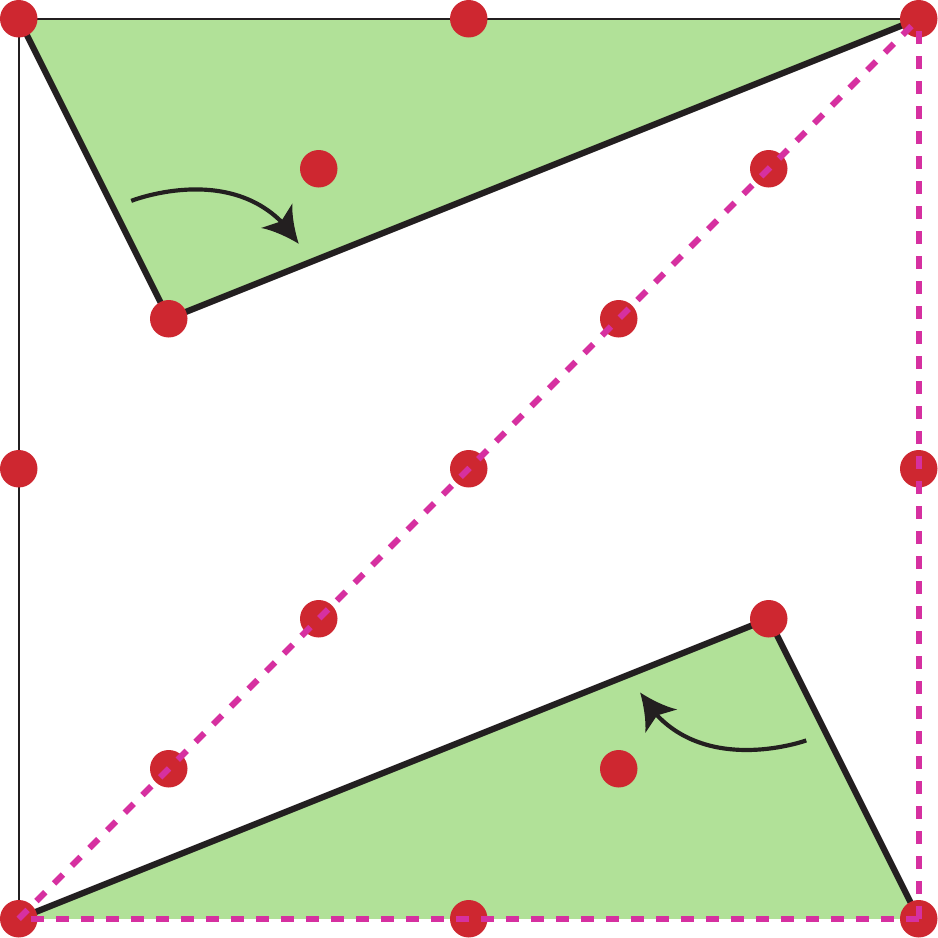}}
\put(-4,0){$O$}
\put(-4,62){$O$}
\put(64,62){$O$}
\put(64,0){$O$}
\put(10,37){$M$}
\put(53,22){$N$}
\put(3,50){$r_1$}
\put(40,50){$r_0$}
\put(57,13){$s_1$}
\put(24,13){$s_0$}
\put(85,0){\includegraphics[width=.42\textwidth]{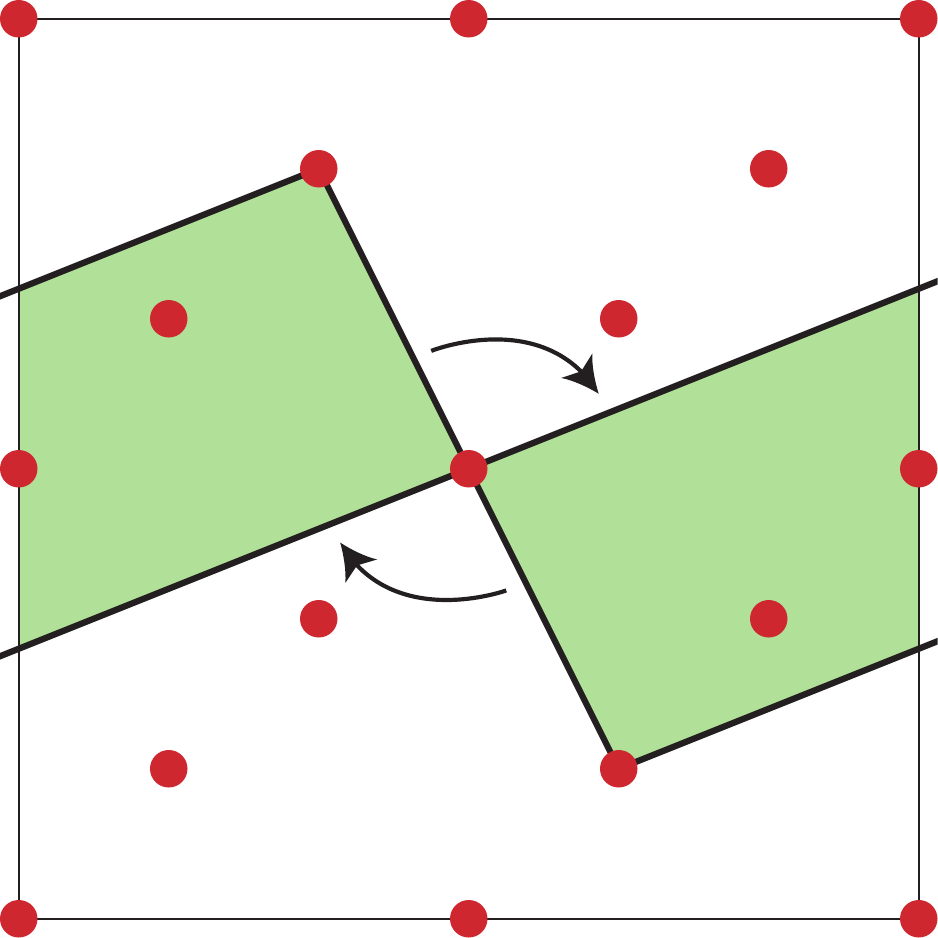}}
\put(118,35){$O$}
\put(123,7){$M$}
\put(106,54){$N$}
\put(123,20){$r_1$}
\put(97,26){$r_0$}
\put(111,43){$s_1$}
\put(134,41){$s_0$}
\end{picture}
\caption{The parallelogram $P$ (green) inside the torus $\TT^2$. The torus is represented as the square $[0,1]^2$ (left) and the square $[-1/2, 1/2]^2$ (right). Red dots are fixed points of $XW$, with $XW = \left(\bsm 3&8\\4&11\esm\right)$ in this example. The action of $XW$ maps $r_1$ to $r_0$ and $s_1$ to $s_0$. If $W$ is not of the form $X^mY^n$ or $Y^nX^m$, the point $N$ lies in the interior of the dotted triangle; otherwise, it lies on its boundary.}
\label{F:Para}
\end{figure}

Let $P$ be the projection of $\widetilde{P}$ onto $\TT^2=\RR^2/\ZZ^2$ (see Figure \eqref{F:Para}). The vertices of $\widetilde{P}$ project onto three points $O, M, N$ which are fixed by the action of $XW$ due to Equation \eqref{eq_minakawa_vector} above. Let $r_1, r_0, s_1, s_0$ denote the projections of the boundary segments of $\widetilde{P}$ onto $\TT^2$. The segments $r_1$ and $r_0$ both connect $O$ to $M$, while $s_1$ and $s_0$ both connect $O$ to $N$. The action of $XW$ maps 
\begin{align*}
    & XW:r_1 \mapsto r_0\\
    & XW:s_1 \mapsto s_0.
\end{align*}

\begin{lemma}\label{lemma_Para}
The interior of $P$ is embedded in $\TT^2$, as well as are embedded the interiors of its sides $r_0, r_1, s_0, s_1$. If $W$ is of the form $XY^n$ or $Y^nX$ for some $n \ge 1$, then $M$ and $N$ coincide in $\TT^2$. Otherwise, the three vertices $O, M$, and $N$ are distinct.
\end{lemma}

\begin{proof}
To show that the interior of $P$ is embedded in $\TT^2$, it suffices to demonstrate that $\widetilde{N}$ lies in the closed triangle bounded by $(0,0), (1, 0), (1,1)$. By the symmetry of the construction, $\widetilde{M}$ then lies in the closed triangle bounded by $(0,0), (1,0), (0,-1)$, ensuring that the triangles $OMO$ and $ONO$ have disjoint interiors in $\TT^2$. 

Recall that all the entries of $W$ are assumed to be positive, and hence the entries of $XW$ are positive integers $a,b,c,d\geq 1$ such that $ad-bc=1$ and $a+d \ge 3$.
In particular, $0 \le a-1 \le a+d-2$, so the first coordinate of $\widetilde{N}$ lies in $[0,1]$. Writing $W = \left(\bsm a'&b'\\c'&d'\esm\right)$ with $a', b', c', d' \geq 1$, we have $XW = (\bsm a'+c'&b'+d'\\c'&d'\esm)$, which implies $0 < c \le a-1$. This confirms that $\widetilde{N}$ lies in the closed triangle defined by $(0,0), (1,0), (1,1)$. 

We now examine when $N$ lies on the boundary of this triangle:
\begin{itemize}
    \item Since $c > 0$, $\widetilde{N}$ cannot lie on the horizontal segment $[(0,0),(1,0)]$. 
    \item If $\widetilde{N}$ lies on the diagonal $[(0,0),(1,1)]$, then $a-1=c$, which implies $a'=1$. Thus $W$ is of the form $Y^n X^m$ for $m,n \ge 1$, and hence $XW = \left(\bsm n+1 & mn+m+1 \\ n & mn+1 \esm\right)$. The trace is $t = mn+n+2$, and $N$ has coordinates $(\frac{1}{m+1}, \frac{1}{m+1})$. The vertices $M$ and $N$ coincide only if $m=1$.
    \item Finally, if $N$ lies on the vertical segment $[(1,0),(1,1)]$ then $a-1=t-2$, which implies $d=1$. This occurs when $XW$ is of the form $X^m Y^n$ for $m \ge 2$ and $n \ge 1$. Here $N$ has coordinates $(1, 1/m)$, and $M=N$ only if $m=2$. 
\end{itemize}
In summary, $P$ fails to be embedded only at $O$ (by construction), and at $N=M$ precisely when $W$ is of the form $XY^n$ or $Y^nX$ for $n \ge 1$.
\end{proof}

\begin{figure}[t]
    \centering
    \includegraphics[width=1\linewidth]{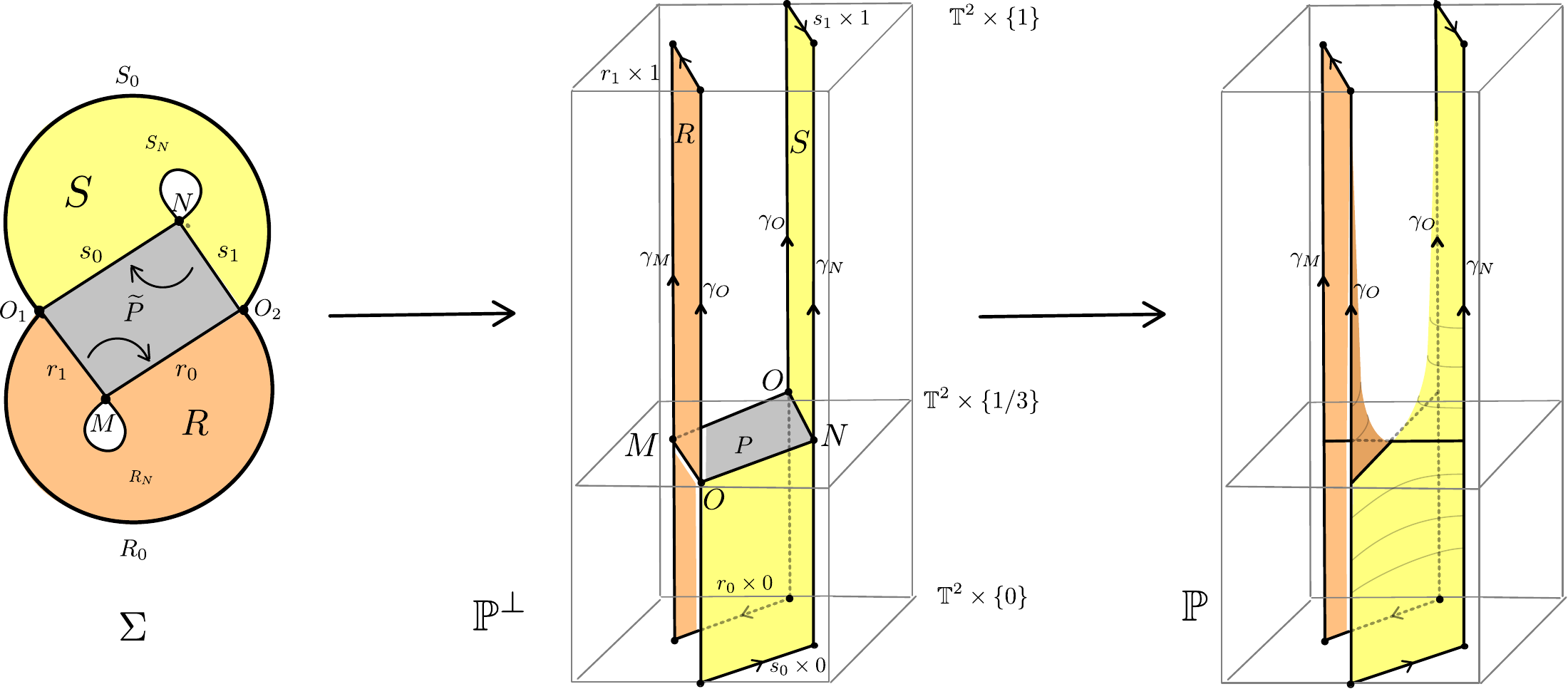}
    \caption{The construction of the immersed pair of pants $\PP^\perp$ (center) and its smoothing~$\PP$ (on the right).}
    \label{F:Pants}
\end{figure}

\paragraph{\textbf{Construction of the pair of pants}}
Recall that $\TT^3_{XW}$ is obtained from $\TT^2 \times [0,1]$ by identifying its boundary components under the relation $(p,1) \sim (XW(p),0)$, and the suspension flow $\Phisus$ is a flow tangent to the segments $p\times[0,1]$. We now construct a pair of pants whose interior is embedded in $\TT^3_{XW}$ and transverse to $\Phisus$. Consider the three periodic orbits of $\Phisus$ given by $\gamma_O := O \times [0,1]$, $\gamma_M := M \times [0,1]$, and $\gamma_N := N \times [0,1]$. Let $\Gamma$ denote the link $\gamma_M \cup \gamma_N \cup \gamma_O$. The boundary of the pair of pants will consist of this link $\Gamma$. 

Consider first the 2-complex $\PP^\perp$ in $\TT^3_{XW}$ obtained as the union of:
\begin{itemize}
    \item the parallelogram $P \times \{1/3\}$ inside the fiber $\TT^2 \times \{1/3\}$,
    \item the vertical strips $r_1 \times [1/3,1]$, $r_0 \times [0,1/3]$, $s_1 \times [1/3,1]$, and $s_0 \times [0,1/3]$.
\end{itemize}
This 2-complex is shown at the center of Figure~\eqref{F:Pants}. Since $s_1 \times \{1\}$ and $r_1 \times \{1\}$ are identified with $s_0 \times \{0\}$ and $r_0 \times \{0\}$ respectively, the complex $\PP^\perp$ can be described as the union of 
and embedded parallelogram $P \times \{1/3\}$ and two embedded strips 
\begin{align*}
 S := (s_1 \times [1/3,1]) \cup (s_0 \times [0,1/3])\\
 R := (r_1 \times [1/3,1]) \cup (r_0 \times [0,1/3])
\end{align*}

Let $\Sigma$ be the abstract orientable surface shown at the left side of Figure~\eqref{F:Pants}, obtained by attaching the strips $S$ and $R$ to the parallelogram $P$ along the sides $s_i$ and $r_i$, $i=0,1$, that is homeomorphic to a pair of pants. The 2-complex $\PP^\perp$ is the image of $\Sigma$ under an immersion, that is not smooth along the attachment segments $r_i \times \{1/3\}$ and $s_i \times \{1/3\}$. Observe that the suspension flow is transverse to $\PP^\perp$ only along the parallelogram and is tangent along the strips, so it fails to be a Birkhoff surface. However, using the technique for smoothing pair of pants described in~\cite{FriedAnosov}, it is possible to make an isotopy of the embedding $\iota:\Sigma\to\TT^3_{XW}$ that fixes the boundary link $\Gamma$, and obtain a smooth Birkhoff surface with boundary on $\Gamma$. 

\begin{lemma}
\label{lemma_Para_2}
The immersion $\iota:\Sigma \to \TT^3_{XW}$ with image $\PP^\perp$ described above can be deformed by isotopies (fixing the link $\Gamma$) into a smooth Birkhoff surface. That is, an immersion that is embedded and transverse to the $\Phisus$-orbits on the interior, and a covering along the boundary. We denote the new immersion with the same symbol $\iota:\Sigma\to \TT^3_{XW}$, and denote by $\PP$ its image. The homology class  $[\iota(\partial\Sigma)]\in\mathrm{H}_1(\TT^3_{XW})$ of its boundary link satisfies 
\begin{align*}
    & [\iota(\partial \Sigma)] = -[\gamma_M] - [\gamma_N] + 2[\gamma_O]\ \text{if}\ W\ \text{is not of the form}\ XY^n\ \text{or}\ Y^nX, \\
    & [\iota(\partial \Sigma)] = -2[\gamma_M] + 2[\gamma_O],\ \text{otherwise}.
\end{align*}
\end{lemma}

\begin{proof}
To check the fact that the 2-complex $\PP^\perp$ can be deformed into a Birkhoff surface, we refer the reader to \cite{FriedAnosov}, where the deformation into Birkhoff surfaces of this of pair of pants is explained. We depicted this construction on the right side of Figure~\eqref{F:Pants}.

To calculate the homology coordinates of the boundary link, let $S_N, S_O$ and $R_M, R_O$ be the boundary arcs of the strips $S$ and $R$ indicated at the left part of Figure \eqref{F:Pants}. The surface $\Sigma$ has three boundary components: $C_N$ corresponding to $S_N$, $C_M$ corresponding to $R_M$, and $C_O$ that corresponds to the union $S_O \cup R_O$. Since the immersion in $M$ of each of the arcs $S_N, R_M, S_O$ or $R_O$ covers its respective periodic orbit ($\gamma_O$, $\gamma_N$ or $\gamma_M$) with degree $\pm1$, we have that $\iota: C_N \to \gamma_N$ and $\iota: C_M \to \gamma_M$ are degree-one coverings, while $\iota: C_O \to \gamma_O$ is a degree-two covering. The orientation of $C_O$ matches the flow orientation on $\gamma_O$, whereas the orientations of $C_N$ and $C_M$ are opposite to the flow. Thus, at the level of homology classes we have that 
\begin{align*}
    & [\iota(C_O)] = 2[\gamma_O]\\
    & [\iota(C_N)] = -[\gamma_N]\\
    & [\iota(C_M)] = -[\gamma_M].
\end{align*}
Finally, using Lemma~\ref{lemma_Para} we can conclude the statement of Lemma \ref{lemma_Para_2} for any matrix $W$.
\end{proof}

\subsection{The Fried sum of $\TT^2_{2/3}$ and $\PP$}\label{S:Cutting}
The torus $\TT^2_{2/3}$ and the pair of pants $\PP$ are two Birkhoff surfaces for the flow $\Phisus$. The former has empty boundary and intersects every orbit (it is a global transverse section), while the latter has non-empty boundary and does not intersect all orbits. These surfaces are transverse to each other and intersect along the two arcs $r_1\times\{2/3\}$ and $s_1\times\{2/3\}$ depicted in Figure\eqref{F:TorusPants}. Their union is not a manifold, but we can construct a new Birkhoff surface from this union using the Fried desingularization mechanism explained at Section~\ref{S:Fried}. We now prove that this section has genus $1$.

\begin{lemma}\label{lemma_Para_3}
The Fried sum $\PP\scup\TT^2_{2/3}$ is the image of a Birkhoff section 
$$\iota:(S,\partial S)\to(\TT^3_{XW},\{\gamma_M,\gamma_O,\gamma_N\})$$
for the suspension flow $(\TT^3_{XW},\Phisus)$, where the surface $S$ has genus 1, and its boundary $\partial S$ satisfy: 
\begin{enumerate}[(a)]
\item If $W$ is not of the form $XY^n$ or $Y^nX$, then $\partial S$ has four components that are mapped onto a set of three different periodic orbits $\{\gamma_M,\gamma_O,\gamma_N\}$ in the following way:
\begin{itemize}
    \item one component of $\partial S$ is mapped onto $\gamma_N$ with multiplicity $-1$,

    \item another component of $\partial S$ is mapped onto $\gamma_M$ with multiplicity $-1$,

    \item the remaining two components are mapped onto $\gamma_O$, each with multiplicity $1$.
\end{itemize}
\vspace*{5pt}
\item If $W$ is of the form $XY^n$ or $Y^nX$, then $\partial S$ has three components that are mapped onto a set of two periodic orbits $\{\gamma_O,\gamma_N=\gamma_M\}$. in the following way
\begin{itemize}
\vspace{5pt}
    \item One component of $\partial S$ is mapped onto $\gamma_M=\gamma_N$ with multiplicity $-2$,

    \item the remaining two components of $\partial S$ are mapped onto $\gamma_O$ with multiplicity $1$. 
\end{itemize}
\end{enumerate}
\end{lemma}

\begin{figure}[t]
\centering
\includegraphics[width=.9\textwidth]{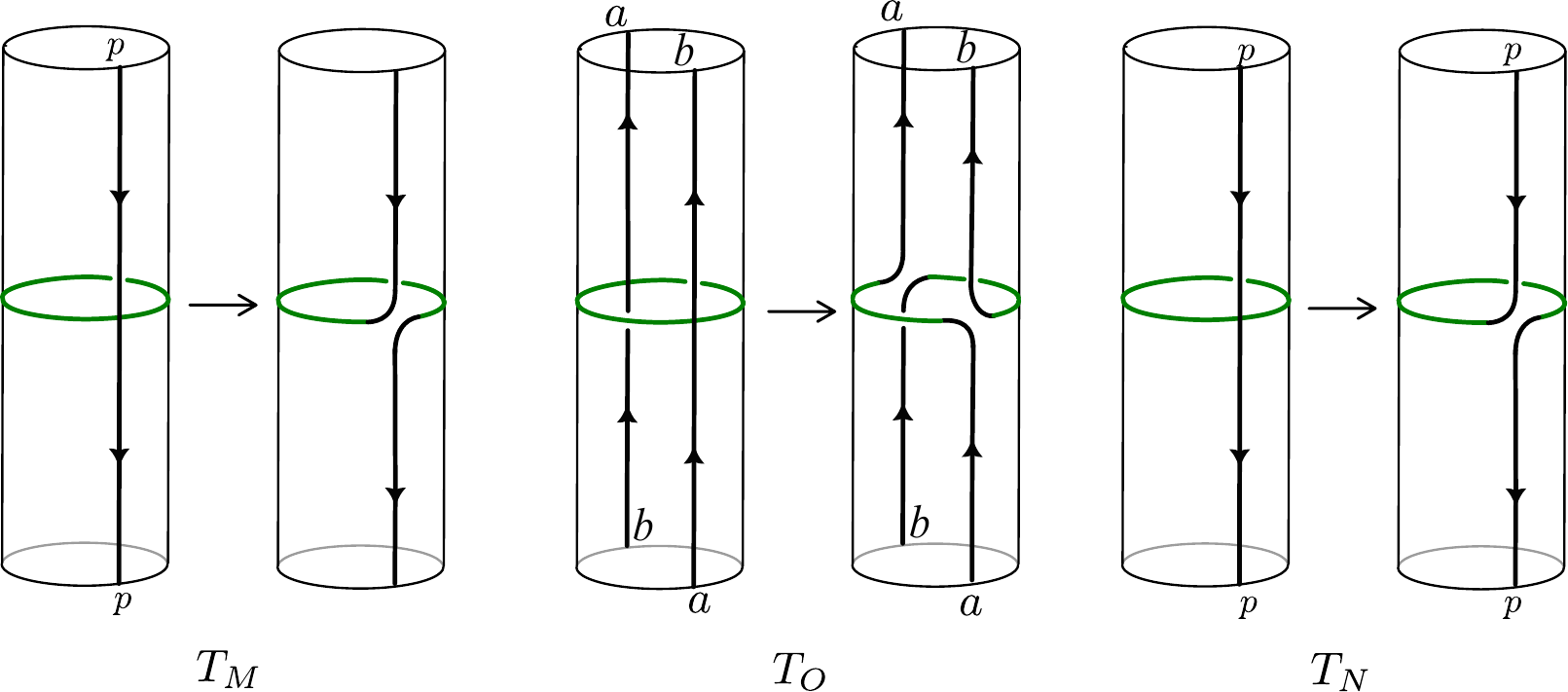}
\caption{Fried desingularization of $\partial\TT^2_{2/3}$ (green) and $\partial\PP$ (black) along on the components $T_M$ (left), $T_O$ (center) and $T_N$ (right) of $\partial\TT^3_{XW, \Gamma}$.}
\label{F:Bord}
\end{figure}

\begin{proof}
First, assume that $W$ is not of the form $XY^n$ or $Y^nX$. 

As explained at Section \ref{S:FriedSum}, the union of $\TT^2_{2/3}$ and $\PP$ can be desingularized, giving rise to a Birkhoff surface $\iota:(S,\partial S)\to(\TT^3_{XW},\Gamma)$ with boundary link $\Gamma=\{\gamma_M,\gamma_O,\gamma_N\}$, whose image we denote by $\PP\scup\TT^2_{2/3}$. We need to determine the genus and boundary components of $S$. 

As noted in Section~\ref{S:Euler}, to compute the Euler characteristic of the embedding and determine its boundary topology, it is convenient to perform the normal blow-up of $\TT^3_{XW}$ along the boundary link $\Gamma$. This is because, if we lift $\PP\scup\TT^2_{2/3}$ to the blown-up manifold, the surface that we obtain is homeomorphic to $S$. Let us denote by $\TT^3_{XW,\Gamma}$ the blow-up manifold, and denote its boundary components by $T_O, T_N,$ and $T_M$, where the subscripts indicate which closed orbit was blown-up.

Let $(\TT^2_{2/3})^*$ and $\PP^*$ be the surfaces lifted to the blow-up manifold $\TT^3_{XW,\Gamma}$. Since $\TT^2_{2/3}$ intersects $\Gamma$ at three interior points and $\PP$ intersects $\Gamma$ only on its boundary, we have that:
$$\chi((\TT^2_{2/3})^*) = -3 \quad \text{and} \quad \chi(\PP^*) = -1.$$
Since the Euler characteristic is additive under the Fried sum of the blown-up Birkhoff sections (cf. Lemma~\ref{lemma_linear_Euler_char}), we obtain:
$$\chi(S)=\chi(\PP\scup\TT^2_{2/3}) = -4.$$ 

Next, we count the boundary components of the embedding $S\to\PP\scup\TT^2_{2/3}$. In the blow-up $\TT^3_{XW,\Gamma}$ the global section $(\TT^2_{2/3})^*$ has three boundary components, each of which is a meridian in $T_N, T_M,$ and $T_O$, depicted in green at Figure\eqref{F:Bord}. The surface $\PP^*$ is a pair of pants, having one boundary component in each boundary torus, depicted in black at Figure\eqref{F:Bord}. 
Observe that:
\begin{itemize}
    \item Since $\PP^*$ has multiplicity $-1$ along $\gamma_N$, the boundary component $\PP^*\cap T_N$ intersects the meridian $(\TT^2_{2/3})^*\cap T_N$ with multiplicity $-1$. The desingularization of the union of these two curves produces a single closed curve in $T_N$ wrapping once in the negative vertical direction (see Figure~\eqref{F:Bord}, right). The same occurs in $T_M$, producing one component wrapping once in the negative vertical direction (Figure~\eqref{F:Bord}, left). 

    \item Since $\PP^*$ has multiplicity $2$ along $\gamma_O$, its boundary component $\PP^*\cap T_O$ intersects the meridian $(\TT^2_{2/3})^*\cap T_O$ with multiplicity $2$. The desingularization of this union results in two closed curves in $T_O$, each wrapping once in the positive vertical direction (Figure~\eqref{F:Bord}, middle). 
\end{itemize}
Consequently, $S\simeq(\TT^2_{2/3})^*\scup\PP^*$ has four boundary components. Given that its Euler characteristic is $-4$, we conclude the surface has genus $1$.

In the case where $W=XY^n$ or $Y^nX$, recall from Lemma~\ref{lemma_Para} that $\gamma_N=\gamma_M$, meaning $\TT^3_{XW,\Gamma}$ has only two boundary components. The Euler characteristic $\chi((\TT^2_{2/3})^*)$ is $-2$, resulting in $\chi(S)=\chi(\PP\scup\TT^2_{2/3}) = -3$. In the blow-up, $\PP^*$ has two boundary components in $T_M$, each intersecting the meridian $(\TT^2_{2/3})^*\cap T_M$ with multiplicity $-1$. Desingularization produces one curve in $T_M$ wrapping twice in the negative vertical direction. The desingularization along $T_O$ remains identical to the previous case. Thus, $S$ has three boundary components and genus $1$.
\end{proof}

\subsection{Computing the first-return map}\label{S:FirstReturn}
We now proceed to compute the first-return map $f_S$ induced by the flow $\Phisus$ on the Birkhoff section $\iota:S\to \PP\scup\TT^2_{2/3}$. 

\begin{remark}[Blow-down of the Birkhoff section]
As explained at the introduction, associated with the Birkhoff section $\iota:S\to \PP\scup\TT^2_{2/3}$, there is a closed surface $\hat{S}$ obtained by collapsing each boundary component of $S$ into a point. The first-return map on $\textrm{int}(S)$ induces a homeomorphism $\hat{f}_S:\hat{S}\to\hat{S}$, where the collapsed boundary components of $S$ become periodic points of $\hat{f}_S$. Since $\Phisus$ is an Anosov flow, its stable and unstable foliations induce a pair of transverse foliations on $\textrm{int}(S)$, that are invariant and uniformly contracted/expanded under the first-return application. These, in turn, extend into a pair of transverse foliations on $\hat{S}$, possibly with \emph{prong-like singularities} at the collapsed boundary points, that are contracted/expanded under the action of $\hat{f}_S$. Under the assumption that the center-stable and center-unstable foliations of the flow are orientable, it is possible to see that the number of prongs at each singularity is greater than one. Therefore, the blow-down application $\hat{f}_S:\hat{S}\to\hat{S}$ is a pseudo-Anosov homeomorphism on a closed surface (see, e.g., \cite{FriedAnosov}).
\end{remark}

In our case, since $\hat{S}$ is a closed surface of genus $1$ (Euler characteristic $\chi(\hat{S})=0$), and since every singularity of the foliation has non-positive index (at least two prongs) then these two foliations have no singularities at all. This implies that $\hat{f}_S$ is conjugated to an Anosov map on the torus $\RR^2/\ZZ^2$, defined by a matrix $A \in \SLZ$ (cf. Lemme 1, page 70 at \cite{GhysGV}). Our goal is to determine this matrix.

\begin{lemma}
\label{lemma_Para_4}
Let $\iota:S\to \PP\scup\TT^2_{2/3}$ be the genus-one Birkhoff section for the suspension Anosov flow $(\TT^3_{XW},\Phisus)$ constructed in Lemma~\ref{lemma_Para_3}. On the blow-down surface $\hat{S}$, the first-return map induces a homeomorphism $\hat{f}_S$ conjugated to the action of the matrix $W$ on $\RR^2/\ZZ^2$. 
\end{lemma}

\begin{figure}[t]
\centering
\includegraphics[width=\textwidth]{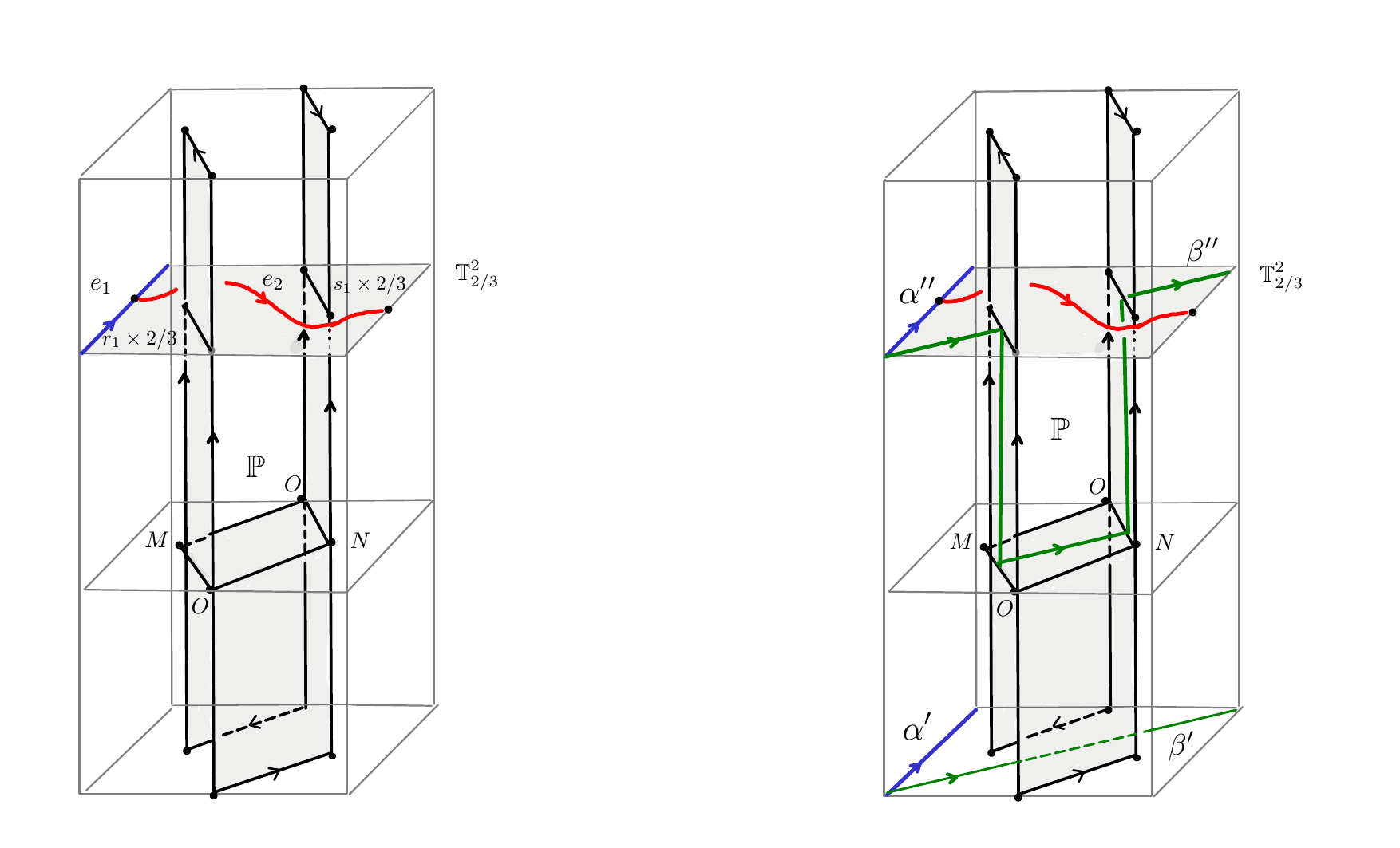}
\caption{The torus $\TT^2_{2/3}$ and the pair of pants $\PP$ intersect along the segments $r_1\times\{2/3\}$ and $s_1\times\{2/3\}$. After applying Fried desingularization, their union gives a genus-one Birkhoff section $\PP\scup\TT^2_{2/3}$ for the flow $\Phisus$ in $\TT^3_{XW}$. The right part illustrates the curves $\alpha'\times\{0\}=(1,0)$ (blue) and $\beta'\times\{0\}=(1,1)$ (green) from Section~\ref{S:FirstReturn}, and their images when pushed toward $\PP\scup\TT^2_{2/3}$ along the flow $\Phisus$.}
\label{F:TorusPants}
\end{figure}

\begin{proof}
Since $\hat{f}_S$ is conjugated to a linear Anosov diffeomorphism on the torus induced by some matrix $A\in\SLZ$, to determine the conjugacy class of $A$ we only need to determine the conjugacy class of the homology action $[\hat{f}_S]:H_1(\hat{S})\to H_1(\hat{S})$, since then $[\hat{f}_S]$ is conjugated to $A$.  

Let us describe first a basis of $\mathrm{H}_1(\hat{S})$ used to express the homology action of $\hat{f}_S$. Consider a pair 
$\{\varepsilon_1,\varepsilon_2\}$ of oriented simple closed curves in $\TT^2$ such that:
\begin{itemize}
    \item $\varepsilon_i$ is disjoint from the union $r_1\cup s_1$ for both $i=1,2$,
    \item $[\varepsilon_1]=(1,0)$ and $[\varepsilon_2]=(0,1)$ in the standard basis of $H_1(\TT^2)$.
\end{itemize}
The choice of these curves is always possible since $r_1\cup s_1$ is contractible, and we may assume that they are transverse and intersect at only one point. The embeddings 
$$e_i = \varepsilon_i \times \{2/3\},\ i=1,2$$ 
yield a pair of simple closed curves in the interior of $\PP\scup\TT^2_{2/3}$ that intersect transversely at a single point (see Figure \eqref{F:TorusPants}). In the blow-down surface $\hat{S}$, the pair $\{e_1, e_2\}$ retains this property, and since $\hat{S}$ has genus equal to $1$, this pair constitutes a basis for $\mathrm{H}_1(\hat{S})$. We make an identification of $\mathrm{H}_1(\hat{S})$ with $\ZZ^2$ by setting:
\begin{itemize}
\item $[e_1]=(1,0)$, $[e_2]=(0,1)$
\end{itemize}

Consider:
\begin{itemize}
    \item a simple closed curve $\alpha$ in $\TT^2_{2/3}$ that avoids $(r_1 \cup s_1) \times \{2/3\}$ and has homology coordinates 
    $$[\alpha]=(d,-c),$$

    \item a simple closed curve $\beta$ in $\TT^2_{2/3}$ that avoids $(r_1 \cup s_1) \times \{2/3\}$ and has homology coordinates 
    $$[\beta]=(d-b,a-c)$$
\end{itemize}
This choice is possible since $r_1 \cup s_1$ is contractible. To calculate the first-return map onto the Birkhoff section, we have to push the curves $\alpha$ and $\beta$ using the forward action of the flow $\Phisus$ until they return to the Birkhoff section, and then express its coordinates in the basis $\{e_1, e_2\}$. We perform this operation in two steps:

\begin{enumerate}

\item First, push $\alpha$ and $\beta$ (both contained in the torus section $\TT^2_{2/3}$) using the forward $\Phisus$ action, until they reach the torus $\TT^2_1$. Once in this torus, this pair of curves is identified with another pair $\alpha' \times \{0\}$ and $\beta' \times \{0\}$ contained in $\TT^2_0$ using the matrix $XW$. Observe that the coordinates $[\alpha]$ and $[\beta]$ in the standard basis of $H_1(\TT^2_1)$ coincide with its coordinates in the basis $\{\varepsilon_1,\varepsilon_2\}$, and hence:
\begin{align*}
&[\alpha'] =XW\cdot [\alpha]= \bpm a & b \\ c & d \epm \bpm d \\ -c \epm = \bpm 1 \\ 0 \epm, \\
&[\beta'] = XW\cdot [\beta] \bpm a & b \\ c & d \epm \bpm d-b \\ a-c \epm = \bpm 1 \\ 1 \epm.
\end{align*}

\item Second, we push the curves $\alpha'\times\{0\}$ and $\beta'\times\{0\}$ contained in the section $\TT^2_0$ using the forward $\Phisus$ action, until they reach $\PP\scup\TT^2_{2/3}$. Observe that
\begin{itemize}
\item When pushing 
$$\Phisus^t:\alpha'\times\{0\}\mapsto\alpha'\times\{t\},$$
the image curves never intersect the pair of pants $\PP$, when 
$0\leq t\leq 1$. Hence, $\alpha'\times\{0\}$ projects onto the curve $\alpha''=\alpha'\times\{2/3\}$ when pushed using the $\Phisus$-action. 
\item When pushing 
$$\Phisus^t:\beta'\times\{0\}\mapsto\beta'\times\{t\},$$
observe that the image curves enter into a \emph{tunnel} formed by $\PP$ for values $t>1/3$ and an arc of the curve gets ``trapped'' below the polygon $P\times\{1/3\}$ (see Figure~\eqref{F:TorusPants}). Therefore, the simple closed curve $\beta''$ in $\PP\scup\TT^2_{2/3}$ obtained by projecting $\beta'\times\{0\}$ using the $\Phisus$-action, does not intersect the curve $(0,1)$ in $\TT^2_{2/3}$ (even though $\beta'$ intersects $(0,1)$ in $\TT^2_0$). 
\end{itemize}
\end{enumerate}

Putting all together, in the homology basis $\{e_1=(1,0), e_2=(0,1)\}$ the first return map takes 
\begin{align*}
    & [\alpha]\mapsto[\alpha'']= \left(\begin{smallmatrix} 1 \\ 0 \end{smallmatrix}\right)\\
    & [\beta]\mapsto[\beta'']=\left(\begin{smallmatrix} 0 \\ 1 \end{smallmatrix}\right).
\end{align*}
Therefore, the homology action $[\hat{f}_S]$ on $H_1(\hat{S})$ of the first-return map on the Birkhoff section is given by:
\[
\bpm d & d-b \\ -c & a-c \epm^{-1} = \bpm a-c & b-d \\ c & d \epm = W.
\]
\end{proof}


\subsection{The case of $YW$}\label{S:case_YW}
Consider again a hyperbolic matrix $W$ with positive entries, and let $(\TT^3_{YW}, \Phisus)$ be the corresponding suspension flow. The procedure for constructing a genus-one Birkhoff section with a first-return map conjugate to $W$ is entirely analogous to the one described in Sections~\ref{S:The_pair_of_pants}, \ref{S:Cutting}, and \ref{S:FirstReturn}, with the following modifications:

\begin{itemize}
    \item The vector $v \in \RR^2$ from Section \ref{S:The_pair_of_pants} is replaced by $v = \left(\begin{smallmatrix} -b/(t-2) \\ (a-1)/(t-2) \end{smallmatrix}\right)$. A direct calculation yields:
    \begin{equation}
    \bpm a & b \\ c & d \epm \bpm -b/(t-2) \\ (a-1)/(t-2) \epm = \bpm -b/(t-2) \\ (1-d)/(t-2) \epm = \bpm -b/(t-2) \\ (a-1)/(t-2) \epm - \bpm 0 \\ 1 \epm.
    \end{equation}
    This implies that $YW \cdot v = v - (0,1)$. We then define a parallelogram $\widetilde{P}_{YW}$ in $\RR^2$ with the following sides and vertices: 
    \begin{align*}
    &\tilde{r}_1 = \{t \cdot v : 0 \le t \le 1\}, &&\widetilde{O}_1 = (0,0),\\
    &\tilde{r}_0 = YW \cdot \tilde{r}_1 + (0,1), &&\widetilde{M} = \left(\frac{-b}{t-2}, \frac{a-1}{t-2}\right),\\
    &\tilde{s}_1 = -\tilde{r}_1 + (0,1), &&\widetilde{O}_2 = (0,1),\\
    &\tilde{s}_0 = -YW \cdot \tilde{r}_1, &&\widetilde{N} = \left(\frac{b}{t-2}, \frac{d-1}{t-2}\right).  
    \end{align*}
    As in Section~\ref{S:The_pair_of_pants}, the action of the matrix $YW$ maps $r_1 \mapsto r_0$ and $s_1 \mapsto s_0$. Geometrically, the polygon $\widetilde{P}_{YW}$ (right part of Figure \eqref{F:Para}) corresponds to a $\pi/2$-rotation of the one shown at the left part of Figure~\eqref{F:Para}.

    \item The construction of the pair of pants and the Fried sum follows the same logic as in Sections~\ref{S:The_pair_of_pants} and \ref{S:Cutting}. The only difference lies in the degenerate cases for Lemmas~\ref{lemma_Para}, \ref{lemma_Para_2}, and \ref{lemma_Para_3}, which here occur when $W = X^n Y$ or $W = Y X^n$ for $n \ge 1$.

    \item The computation of the first-return map follows the steps in Section~\ref{S:FirstReturn}. In the $XW$ case, adding the pair of pants $\PP$ to the global section $\TT^2_{2/3}$ effectively "subtracts" the matrix $X$ from the left. In the present case, the corresponding pair of pants $\PP$ for $(\TT^3_{YW}, \Phisus)$ is rotated by $\pi/2$ relative to the previous construction, and its addition effectively "subtracts" the matrix $Y$ from the left of the first-return map.
\end{itemize}


\subsection{Proof of the main results}
\label{S:proof_thmB}

We begin with the proof of Theorem~\ref{T:DistanceBis}. Given a hyperbolic matrix $W\in\SLZ$ with positive entries, Lemmas~\ref{lemma_Para_3} and \ref{lemma_Para_4} (including their modifications in Section~\ref{S:case_YW}) establish that the flow $(\TT^3_{ZW}, \Phisus)$, where $Z \in \{X, Y\}$ as before, admits a genus-one Birkhoff section. This section has a boundary link supported by at most three periodic orbits of $\Phisus$, and its first-return map is conjugated to $W$ in the blow-down surface. This concludes the proof of Theorem~\ref{T:DistanceBis}. 

Theorem~\ref{T:Distance} follows as a direct corollary of Theorem~\ref{T:DistanceBis}. 

To prove Theorem~\ref{T:Trace}, consider a matrix $A$ in $\SLZ$ with $\tr(A) > 3$. Since the matrix $A$ is hyperbolic, as explained at the introduction, it is conjugated to a matrix $W$ that is a finite word in $\{X=(\begin{smallmatrix}1&1\\0&1\end{smallmatrix}), Y=(\begin{smallmatrix}1&0\\1&1\end{smallmatrix})\}$ containing both letters, that is unique up to cyclic permutations. It is enough to prove the statement of Theorem \ref{T:Trace} for the suspension flow $(\TT^3_{W},\Phisus)$, since this last one is conjugated to $(\TT^3_{A},\Phisus)$. Since $\tr(A)=\tr(W)$ is strictly bigger than $3$, then the word $W$ contains at least three letters in total and it can be written in the form $W=ZB$, where $Z \in \{X, Y\}$ and $B$ is a smaller word containing both letters (hence hyperbolic). Applying Theorem~\ref{T:DistanceBis} provides a genus-one Birkhoff section for $(\TT^3_{ZB}, \Phisus)$ with a first-return map given by $B$, and hence an almost equivalence between $(\TT^3_W,\Phisus)$ and $(\TT^3_B,\Phisus)$. Since $B$ is a hyperbolic matrix with positive entries (it equals a positive word in $\{X,Y\}$) then by simple calculation we see that $3 \le \tr(B) < \tr(ZB) = \tr(A)$. This proves Theorem~\ref{T:Trace}.

Finally, Theorem~\ref{T:Main} is obtained by an iterative application of Theorem~\ref{T:Trace}. Given any matrix $A$ with $\tr(A) > 3$, Theorem \ref{T:Trace} provides an almost equivalence 
\begin{equation*}
    (\TT^3_A,\Phisus)\rightarrow (\TT^3_B,\Phisus)
\end{equation*}
where the matrix $B$ has trace satisfying $3\leq \tr(B)<\tr(A)$. By defining $A_0=A$, $A_1=B$ and iterating the process, we obtain a sequence of almost-equivalences
\begin{equation*}
    (\TT^3_{A_{i-1}},\Phisus)\rightarrow (\TT^3_{A_{i}},\Phisus),\ i\geq 1
\end{equation*}
satisfying $3\leq \tr(A_i)<\tr(A_{i-1})$ that must necessarily stop in a finite number of steps when $\tr(A_i)=3$, in which case $A_i=XY=(\begin{smallmatrix}2&1\\1&1\end{smallmatrix})$. By concatenating all the almost-equivalences in the process, we obtain an almost equivalence between $(\TT^3_A,\Phisus)$ and $(\TT^3_{XY},\Phisus)$. This concludes the proof of Theorem \ref{T:Main}.


\begin{proof}[Proof of Corollary~\ref{C:Distance}]
Let $A$ be a hyperbolic matrix with positive trace, and let $W$ be a finite word in $\{X,Y\}$ conjugated to $A$. Let $m=m(A)=m(W)$ be the word length. We will make the proof using $(\TT^3_W,\Phisus)$, which is orbit equivalent to $(\TT^3_A,\Phisus)$. 

Write $W=ZW_1$, where either $Z=X$ or $Z=Y$, and $W_1$ is a smaller word in $\{X,Y\}$ containing both letters. Since $W_1$ has positive entries, then Theorem~\ref{T:DistanceBis} apply and the edges connecting the suspensions $(\TT^3_{ZW_1},\Phisus)$ and $(\TT^3_{W_1},\Phisus)$ are obtained by removing a set of three periodic orbits on each flow (cf. Lemma \ref{lemma_Para_3}-(a)). 
Let us denote by $\{\gamma_O,\gamma_M,\gamma_N\}$ the set of orbits removed from the first flow, and $\{\delta_O,\delta_M,\delta_N\}$ the one removed from the second. The almost-equivalence obtained obtained after removing these sets of periodic induces a bijection 
$$\{\gamma_O,\gamma_M,\gamma_N\}\to\{\delta_O,\delta_M,\delta_N\},$$ 
so we have chosen the subindices to respect this identification.   
Along with the construction, the orbit $\gamma_O$ is defined to be the periodic orbit of $(\TT^3_{ZW_1},\Phisus)$ associated to the fixed point $(0,0)$ of the matrix $ZW_1$ acting on $\TT^2$, and $\gamma_N$, $\gamma_M$ correspond to other two fixed points. With respect to the other flow, the orbit $\delta_O$ corresponds to a periodic point of period $2$ for the action of $W_1$ on $\TT^2$, while $\delta_N$ and $\delta_M$ correspond to a pair of fixed points. 

If we want to connect now $(\TT^3_{W_1},\Phisus)$ with some $(\TT^3_{W_2},\Phisus)$ by deleting a letter on the left of the word $W_1$, observe that after applying a self-orbit equivalence, we may assume that the orbit $\delta_N$ (or $\delta_M$) corresponds to the fixed point $(0,0)$ of $W_1$ acting on $\TT^2$, and hence we apply the construction of Theorem \ref{T:DistanceBis} again. Observe that, along the path 
$$(\TT^3_{W},\Phisus)\to (\TT^3_{W_1},\Phisus)\to (\TT^3_{W_2},\Phisus)$$
we need to remove five orbits from $(\TT^3_{W},\Phisus)$ to connect it with $(\TT^3_{W_2},\Phisus)$. Actually, along the first connection we have to remove three periodic orbits $\gamma_O,\gamma_N,\gamma_M$; while in the second, one of the three removed curves is the image $\gamma_N\mapsto\delta_N$ under the first almost-equivalence, adding only two curves in $(\TT^3_{W},\Phisus)$ to be removed. This is completely general: If we use the procedure of Theorem \ref{T:DistanceBis} to connect two words by successively removing one letter on the left of the first one, we have to remove three periodic orbits due to the first letter deletion, and two more for each other letter deletion.  

Finally, observe that every word can be reduced to either~$X^2Y$ or $
XY^2$ by deleting $m-3$ letters, and each of the latter is connected to $XY$ by removing only two orbits (cf. Lemma \ref{lemma_Para_3}-(b)), where one of them corresponds to the fixed point $(0,0)$ of the action on $\TT^2$. Hence, when connecting $(\TT^3_{W}\Phisus)$ to $(\TT^3_{XY},\Phisus)$ by deleting letters on the left, the last step only add one curve to be removed from $(\TT^3_{W}\Phisus)$.
Putting all together, we obtain the bound
$$\dgood\left((\TT^3_W,\Phisus),(\TT^3_{XY},\Phisus)\right)\leq 2m-4.$$

For the second bound involving $\tr(A)=\tr(W)$, observe that in the path of almost equivalences described above, since each $W_i$ satisfies $3\leq\tr(W_i)<\tr(W_{i-1})$, then there can be at most $\tr(W)-3$ steps in the path connecting $W$ with $XY$. The first step contributes with three periodic orbits to remove, while the remaining $\tr(W)-4$ contributes with at most two curves. Hence, we obtain
$$\dgood\left((\TT^3_W,\Phisus),(\TT^3_{XY},\Phisus)\right)\leq 2\tr(A)-5.$$
\end{proof}

\medskip
We conclude with an instructive remark regarding the original argument by Minakawa for proving Theorem~\ref{T:Trace}:

\begin{remark}
Given a matrix $A$ with $\tr(A) > 3$, it can always be written as $A = ZB$ where $Z \in \{X, Y\}$ and $B$ has positive entries, and by the construction at Lemma \ref{lemma_Para_3} the flow $(\TT^3_A, \Phisus)$ admits a genus-1 Birkhoff section with a first-return map conjugate to $B$. The central point is to show that $B$ has \textbf{fewer fixed points} than $A$. By the formula for the number of fixed points of an Anosov automorphism, $|\mathrm{Fix}(E)| = |\tr(E)-2|$ for $E \in \SLZ$, the inequality $|\mathrm{Fix}(B)| < |\mathrm{Fix}(A)|$ implies $\tr(B) < \tr(A)$.

To see that $B$ has \textbf{fewer fixed points} than $A$, note that the surface $\TT^2_{2/3} \scup \PP$ intersects every closed orbit of the suspension flow at least as many times as $\TT^2_{2/3}$. Consequently, the number of periodic orbits intersected exactly once by $\TT^2_{2/3} \scup \PP$ is less than or equal to those intersected by $\TT^2_{2/3}$. The equality cannot hold: while the orbit $\gamma_O$ corresponds to a fixed point for the first-return map $A$ on $\TT^2_{2/3}$, the two boundary components of the Birkhoff section arriving at $\gamma_O$ produce an orbit of order 2 for $B$ when collapsed into points.\footnote{As shown in Figure~\eqref{F:Para}, the number of fixed points of $B$ is generally much smaller than that of $A$, since any fixed point for $A$ located inside the interior of $P_A$ becomes a periodic point of higher period for $B$.} This completes the proof of Minakawa's Theorem~\ref{T:Trace}.
\end{remark}


\bibliographystyle{siam}

\begin{thebibliography}{15}



\bibitem[BBY17]{BBY}{\sc B\'eguin} Fran\c cois, {\sc Bonatti} Christian, {\sc Yu} Bin : \href{https://doi.org/10.2140/gt.2017.21.1837}{Building Anosov flows on 3-manifolds}, {\it Geom. Topol.} {\bf 21} (2017), 1837--1930. 

\bibitem[BaF22]{BF}{\sc Barbot} Thierry \& {\sc Fenley} Sergio: \href{https://arxiv.org/abs/2205.02495}{Orbital equivalence classes of finite coverings of geodesic flows}, {\it preprint} (2022), arXiv:2205.02495.

\bibitem[Bir1917]{Birkhoff}{\sc Birkhoff} George D : \href{https://doi.org/10.1073/pnas.3.4.314}{Dynamical systems with two degrees of freedom}, {\it Trans.\ Amer.\ Math.\ Soc.\ }{\bf 18} (1917), 199--300. 

\bibitem[BoI23]{Bonatti-Iako}{\sc Bonatti} Christian and {\sc Iakovoglou} Ioannis : \href{https://doi.org/10.1017/etds.2021.170}{Anosov flows on $3$ -manifolds: the surgeries and the foliations}, {\textit Ergod. Th. \& Dynam. Sys.} {\bf 43} (2023), n4, 1129--1188.

\bibitem[Bru94]{Brunella}{\sc Brunella} Marco : \href{http://eudml.org/doc/73338}{On the discrete Godbillon-Vey invariant and Dehn surgery on geodesic flow}, {\it Ann.\ Fac.\ Sc.\ Toulouse S\'er.\ 6} {\bf 3} (1994), 335--344.

\bibitem[Deh11]{Lorenz}{\sc Dehornoy} Pierre : \href{https://doi.org/10.4171/LEM/57-3-1}{Les n\oe uds de Lorenz}, {\it Enseign. Math.} {\bf 57} (2011), 211--280. 

\bibitem[Deh13]{Commensurability}{\sc Dehornoy} Pierre : \href{https://doi.org/10.1016/j.crma.2013.02.012}{Almost-commensurability of 3-dimensional Anosov flows}, \emph{C. R. Math.} {\bf 351} (2013), 127--129. 

\bibitem[Deh15]{GenusOne}{\sc Dehornoy} Pierre : \href{https://doi.org/10.1017/etds.2014.14}{Genus one Birkhoff sections for geodesic flows}, \emph{Ergod. Th. \& Dynam. Sys.} {\bf 35} (2015), 1795--1813. 

\bibitem[Fen94]{Fenley_Anosov-flows}{\sc Fenley} Sergio : \href{https://doi.org/10.2307/2946628}{Anosov Flows in 3-Manifolds}, {\it Ann. of Math. (2)} {\bf 139} (1994), n1, 79--115.

\bibitem[FrW80]{FW}{\sc Franks} John, {\sc Williams} Robert : \href{https://doi.org/10.1007/BFb0086986}{Anomalous Anosov flows}, in {\it Global theory of Dyn. Systems}, Lecture Notes in Math. {\bf 819} (1980) Springer.

\bibitem[Fri83]{FriedAnosov}{\sc Fried} David : \href{https://doi.org/10.1016/0040-9383(83)90015-0}{Transitive Anosov flows and pseudo-Anosov maps}, {\it Topology} {\bf 22} (1983), 299--303.

\bibitem[Ghy84]{GhysAnosov}{\sc Ghys} \'Etienne : \href{ https://doi.org/10.1017/S0143385700002273}{Flots d'Anosov sur les 3-variétés fibrées en cercles}, {\it Ergod. Th. \& Dynam. Sys.} {\bf 4} (1984), n1, 67--80. 

\bibitem[Ghy87]{GhysGV}{\sc Ghys} \'Etienne : \href{https://doi.org/10.5802/aif.1111}{Sur l'invariance topologique de la classe de Godbillon-Vey}, {\it Ann.\ Inst.\ Fourier} {\bf 37} (1987), 59--76. 

\bibitem[Goo81]{Goodman_surgery}{\sc Goodman} Sue : \href{https://link.springer.com/chapter/10.1007/BFb0061421}{Dehn surgery on Anosov flows}, in {\it Geometric dynamics (Rio de Janeiro, 1981)}, 300--307, Lecture Notes in Math., 1007, Springer, Berlin

\bibitem[Gro76]{Gromov}{\sc Gromov} Mikha\"il : {Three remarks on geodesic dynamics and fundamental group}, preprint SUNY (1976), reprinted in {\it Enseign.\ Math.\ (2)} {\bf 46 } (2000), 391--402.
 

\bibitem[HaT80]{HandelThurston}{\sc Handel} Mikael, {\sc Thurston} William : \href{http://eudml.org/doc/142732}{Anosov flows on new three manifolds}, {\it Invent. Math.} {\bf 59} (1980), 95--103.


\bibitem[Has90]{Hashiguchi_first-return}{\sc Hashiguchi} Norikazu : On the Anosov diffeomorphisms corresponding to geodesic flows on negatively curved closed surfaces, {\it J. Fac. Sci. Univ.} Tokyo, {\bf 37} (1990), 485-494.

\bibitem[HaM13]{HM}{\sc Hashiguchi} Norikazu and {\sc Minakawa} Hiroyuki : \href{https://doi.org/10.2969/aspm/07210367}{Genus One Birkhoff Sections for the Geodesic Flows of Hyperbolic 2-Orbifolds}, in {\it Steve Hurder, Takashi Tsuboi, Taro Asuke, Shigenori Matsumoto, and Yoshihiko Mitsumatsu. Geometry, Dynamics, and Foliations 2013 : In Honor of Steven Hurder and Takashi Tsuboi on the Occasion of Their 60th Birthdays}, 2017.

\bibitem[KaH95]{Katok-Hasselblatt}{\sc Katok} Anatole and {\sc Hasselblatt} Boris : \textit{Introduction to the Modern Theory of Dynamical Systems}, Encyclopedia of Mathematics and its Applications, 54, Cambridge Univ. Press, Cambridge, 1995

\bibitem[Kir97]{Kirby}{\sc Kirby} Rob : \href{https://people.brandeis.edu/~ruberman/K2.pdf}{\it Problems in low-dimensional topology}, {Geometric Topology Monographs}, vol. 2, International Press, 1997, pp. 35--473.

\bibitem[Mar24]{Marty}{\sc Marty} Théo : \href{https://doi.org/10.4171/JEMS/1619}{Skewed Anosov ﬂows in dimension 3 are Reeb-like}, {\it J. Eur. Math. Soc.} (2025), published online first.

\bibitem[Min13]{Minakawa}{\sc Minakawa} Hiroyuki : \href{https://www.ms.u-tokyo.ac.jp/video/conference/2013GF/cf2013-038.html}{Genus one Birkhoff sections for suspension Anosov flows}, {\it talk at the conference ``Geometry, Dynamics, and Foliations 2013''}, Tokyo. 

\bibitem[Pau25]{Paulet}{\sc Paulet} Neige : \href{https://doi.org/10.3934/jmd.2025002}{Anosov flows in dimension 3 from gluing building blocks with quasi-transverse boundary}, {\it J. Mod. Dyn.} {\bf 21} (2025), 21--240.

\bibitem[Plan81]{Plante}{\sc Plante} John F : \href{https://doi.org/10.1112/jlms/s2-23.2.359}{Anosov flows, transversely affine foliations, and a conjecture of Verjovsky}, {\it J. London Math. Soc. (2)} {\bf 23} (1981), n2, 359--362.
 
\bibitem[Sha21]{Shannon}{\sc Shannon} Mario : \href{https://arxiv.org/abs/2108.12000}{Hyperbolic models of transitive topologically Anosov flows in dimension three
}, {\it preprint} arXiv:2108.12000 (2021).

\bibitem[Sha24]{Sha24}{\sc Shannon} Mario : \href{https://arxiv.org/abs/2410.23551}{Infinitely-many paths in the graph of Anosov flow}, {\it preprint} arXiv:2410.23551 (2024).

\bibitem[Tom68]{Tomter}{\sc Tomter} Per : \href{https://www.proquest.com/openview/0a0cb284ed9975feaccc2a4198d2cd82/1?pq-origsite=gscholar&cbl=18750&diss=y}{Anosov flows on infra-homogeneous spaces}, in {\it Global Analysis }(Proc. Sympos. Pure Math., Vol.{XIV}, Berkeley, Calif., 1968).

\bibitem[Tsa24a]{Tsang1}{\sc Tsang} Chi Cheuk : \href{https://doi.org/10.1017/etds.2023.105}{Constructing Birkhoff sections for pseudo-Anosov flows with controlled complexity}, {\it Ergod. Th. \& Dynam. Sys.} 44, (2024), n8, 2308--2360.

\bibitem[Tsa24b]{Tsang2}{\sc Tsang} Chi Cheuk : \href{https://doi.org/10.1017/etds.2023.105}{Examples of Anosov flows with genus one Birkhoff sections}, {\it preprint} arXiv:2402.00229 (2024). 

\end{thebibliography}

\end{document}